\documentclass[a4paper,11pt]{article}
\usepackage{epsfig}
\usepackage{amsmath}
\usepackage{amssymb}
\usepackage{amsthm}
\usepackage{bbold}
\usepackage{algorithm2e}
\usepackage{color,xcolor}
\usepackage{tensor}
\usepackage{float}
\usepackage[
	bookmarks,
	colorlinks=true, 
	pdftitle={Hamiltonian with dissipation},
	urlcolor= blue,
	linkcolor= blue,
	citecolor= blue
]{hyperref}
\usepackage[
	style=numeric,
	citestyle=numeric,sorting=nyvt,sortcites=false,
	giveninits=true,
	maxnames=10,
	url=false,eprint=true,doi=false,isbn=false,
	backref=false,
    backend=biber,
	natbib=true,
]{biblatex}
\AtEveryBibitem{\clearfield{note}} 
\AtEveryBibitem{\ifentrytype{book{\clearfield{pages}}{}}}

\usepackage[hyphenbreaks]{breakurl}
\usepackage{cleveref}%

\newcommand{\norm}[1]{\ensuremath{\lVert #1\rVert}}
\newcommand{\Ccal}{\mathcal{C}}
\newcommand{\Ecal}{\mathcal{E}}
\newcommand{\Fcal}{\mathcal{F}}
\newcommand{\Hcal}{\mathcal{H}}
\newcommand{\Jcal}{\mathcal{J}}
\newcommand{\Kcal}{\mathcal{K}}
\newcommand{\Lcal}{\mathcal{L}}

\newcommand{\Pcal}{\mathcal{P}}
\newcommand{\Scal}{\mathcal{S}}
\newcommand{\Tcal}{\mathcal{T}}
\newcommand{\Ucal}{\mathcal{U}}
\newcommand{\bra}[2]{\ensuremath{\left\{#1,#2\right\}}}
\newcommand{\inp}[1]{\left<{#1}\right>}
\newcommand{\inpo}[1]{\left<{#1}\right>_0}
\newcommand{\exd}{\mathsf{d}} 
\newcommand{\contr}{\mathsf{i}} 
\newcommand{\Lie}[1]{\ensuremath{\mathsf{L}_{{#1}}}} 
\newcommand{\dH}{\delta \mathcal{H}}
\newcommand{\dS}{\delta \mathcal{S}}
\newcommand{\dt}{\Delta t}
\newcommand{\Jc}{\mathcal{J}}
\newcommand{\Gc}{\mathcal{G}}
\newcommand{\Jd}{J}
\newcommand{\Gd}{G}
\newcommand{\Jm}{\widehat{J}}
\newcommand{\Gm}{\widehat{G}}

\DeclareMathOperator{\DIV}{div}
\DeclareMathOperator{\curl}{curl}

\DeclareMathOperator{\per}{per}

\numberwithin{equation}{section}

\newtheorem{theorem}{Theorem}[section]
\newtheorem{lemma}[theorem]{Lemma}
\newtheorem{proposition}[theorem]{Proposition}
\theoremstyle{definition}
\newtheorem{definition}{Definition}[section]
\newtheorem{example}{Example}[section]
\theoremstyle{remark}

\newcommand{\email}[1]{\protect\href{mailto:#1}{#1}}

\title{Conformal variational discretisation of infinite dimensional Hamiltonian systems with gradient flow dissipation}

\author{Damiano Lombardi\thanks{Sorbonne Université, Inria and CNRS, UMR 7598 Laboratoire Jacques-Louis Lions, Paris, France. (\email{damiano.lombardi@inria.fr}).}
\and Cecilia Pagliantini\thanks{Dipartimento di Matematica,
			  Universit\`a di Pisa, Pisa,
			  Italy.
  (\email{cecilia.pagliantini@unipi.it}).\\
  C. Pagliantini acknowledges the MIUR Excellence Department Project awarded to the Department of Mathematics, University of Pisa, CUP I57G22000700001, and the INDAM/GNCS 2024 project CUP E53C23001670001.}
}
\date{}

\begin{document}

\maketitle


\begin{abstract}
    Nonconservative evolution problems describe irreversible processes and dissipative effects in a broad variety of phenomena. Such problems are often characterised by a conservative part, which can be modelled as a Hamiltonian term, and a nonconservative part, in the form of gradient flow dissipation. Traditional numerical approximations of this class of problem typically fail to retain the separation into conservative and nonconservative parts hence leading to unphysical solutions. In this work we propose a mixed variational method that gives a semi-discrete problem with the same geometric
    structure as the infinite-dimensional problem. As a consequence the conservation laws and the dissipative terms are retained. Convergence results on the solution are established.
    Numerical tests of the Korteweg-de Vries equation and of the two-dimensional Navier-Stokes equations on the torus and on the sphere are presented to corroborate the theoretical findings.
\end{abstract}

\section{Introduction}


We consider time-dependent partial differential equations characterized by a conservative part, expressed through a Hamiltonian term, and a gradient flow dissipation.
The problem is as follows: we look for $u$ in a Hilbert space $V$ 
such that
\begin{equation}\label{eq:pbm}
\partial_t u = \Jc(u)\dH(u) + \Gc(u)\dS(u),\qquad 
t>0
\end{equation}
supplied by the initial condition $u(0)=u_0\in V$
and suitable boundary conditions encoded in $V$.
In the evolution equation \eqref{eq:pbm},
$\Jc(u)$ characterizes the conservative term and is a Poisson operator, for any $u\in V$, while $\Gc(u)$ characterizes the dissipation and is a positive/negative semi-definite self-adjoint operator. The function $\mathcal{H}: V\rightarrow \mathbb{R}$ plays the role of Hamiltonian, while $\mathcal{S}: V\rightarrow \mathbb{R}$ is a produced/dissipated quantity. We postpone to \Cref{sec:HamGrad} a detailed description of the model problem \eqref{eq:pbm}.

Systems of the form \eqref{eq:pbm} describe irreversible processes and dissipative effects in a broad variety of phenomena.
Examples are dissipative systems characterised by thermal conduction, electric resistivity, viscous dissipation, and phenomena
modelled by gradient flows \cite{jordan1998, mielke23}.
Problems \eqref{eq:pbm} can also model phenomena where energy is conserved but entropy production and viscous heating are allowed. These processes are typical of a large variety of time-dependent systems out of thermodynamic equilibrium such as complex fluids \cite{GEN97b}, geophysical fluids \cite{EGB20}, climate models \cite{VL20}, port-based networks \cite{vdS06}, relativistic hydrodynamics and cosmology models \cite{Ott98}, etc.. The latter class comprises also dissipative Hamiltonian flows and evolution problems ensuing from the Navier-Stokes equations, and kinetic models, like the Boltzmann equation and the Vlasov equation with collisions \cite{Kau84}, where energy, mass, and momentum are conserved and entropy is produced.
%
The formulation \eqref{eq:pbm} includes also, as limit cases, Hamiltonian PDEs, obtained when $\Gc(u)$ vanishes, and gradient flows, when $\Jc(u)$ is zero.

In the numerical discretisation of problem \eqref{eq:pbm} it is of crucial importance to retain the separation into a conservative part and a dissipative part to ensure a physically consistent approximation of the solution.
It is by now well established that formulating numerical approximations compatible with the geometric structure underlying the model problem and not just approximating it yields numerical discretisations with superior stability properties.
The geometric discretisation of canonical Hamiltonian PDEs has received considerable attention, see \cite{LeRe05} for a review and references, and, similarly, structure-preserving approximations of dissipative evolution problems have been proposed both for finite dimensional \cite{QMacL08} and infinite dimensional problems \cite{Egger19}.
By contrast the numerical approximation of noncanonical Hamiltonian systems with gradient flow dissipation is, to the best of our knowledge, largely unexplored.

The bulk of geometric discretisations for problems of the form \eqref{eq:pbm} has focused on specific examples such as the Navier-Stokes equations \cite{Reb07,Palha24}
or a specific subclass, such as port-Hamiltonian systems \cite{Golo04}.
A recent contribution to this topic is the finite element in time discretisation proposed in \cite{AF24} where, although the geometric structure of the problem is not preserved, dissipation laws and conservation of general invariants are enforced via a particular choices of test functions.

In this work we derive a variational method for the numerical \emph{spatial} approximation of evolution equations of the form \eqref{eq:pbm} resulting in a semi-discrete problem with the same geometric structure as the continuous problem, in particular the separation into a conservative and a dissipative part.
The approach we propose in this work is to first introduce a mixed formulation that, under natural regularity assumptions on the differential operators, is equivalent to the classical variational formulation of \eqref{eq:pbm}. 
Then a semi-discrete problem is obtained from a conformal spatial discretisation on a finite dimensional space $V_N\subset V$ of the aforementioned mixed variational formulation. The resulting approximation is shown to be composed of a conservative part characterized by a finite dimensional Poisson operator, and a dissipative part associated to a finite dimensional symmetric positive/negative semi-definite operator. The mixed formulation, at the expense of having two extra unknowns, enables a structure-preserving discretisation with classical low-order finite elements.
Indeed, we prove that: (i) the discretisation of the invariants of motion of \eqref{eq:pbm} yields invariants of motion of the semi-discrete problem; (ii) the kernel of the operators $\Jcal$ contained in the approximation space $V_N$ is preserved; (iii) the equilibria belonging to $V_N$ are equilibria of the semi-discrete problem; and (iv) the semi-discretisation of a metriplectic system, as defined in \Cref{sec:metr}, remains metriplectic.
In terms of accuracy of the approximation, we establish convergence of the solution of the semi-discrete system to the solution of the continuous mixed formulation.

The remainder of the paper is organised as follows. 
In \Cref{sec:HamGrad} we introduce the mathematical framework of evolution equations of the form \eqref{eq:pbm} by first characterizing the properties of Hamiltonian PDEs and, then, of gradient flows.
\Cref{sec:weakform} pertains to the derivation of the variational formulation of problem \eqref{eq:pbm} and its mixed variational formulation.
\Cref{sec:semidisc} is the core part of the paper where we introduce the proposed conformal variational discretisation based on the mixed formulation of the problem, we prove the preservation of the geometric structure of the problem and establish convergence results.
In \Cref{sec:time} we discuss the temporal discretisation of the semi-discrete problem via established symplectic methods and discrete gradients. Several numerical tests on the Korteweg-de Vries equation and on the 2D Navier-Stokes equation are presented in \Cref{sec:num}.
Some concluding remarks are drawn in \Cref{sec:conclusions}.


\section{Hamiltonian systems with gradient flow dissipation}\label{sec:HamGrad}

\subsection{Infinite dimensional Hamiltonian dynamics}\label{sec:Ham}


The definition of Hamiltonian partial differential equations require the introduction of function spaces endowed with a symplectic manifold structure and a characterization of the velocity field describing the Hamiltonian flow.
To this end we recall some important concepts and refer to, e.g. \cite{kuksin00}, for a detailed introduction to the topic.
First we introduce the concept of symplectic Hilbert scales.
\begin{definition}[Hilbert scale]
Let $X_0$ be a separable real Hilbert space with scalar product $\inpo{\cdot,\cdot}$ and a Hilbert basis $\{\psi_k\}_{k\in\mathbb{Z}}$.
Let $\{\theta_k\}_{k\in\mathbb{Z}}$ be a positive sequence such that
$\theta_k\to\infty$ as $|k|\to\infty$.
For any $s\in\mathbb{Z}$, let $X_s$ be a Hilbert space with
norm and scalar product
denoted by $\norm{\cdot}_s$ and $\inp{\cdot,\cdot}_s$, respectively, and
defined as
\begin{equation*}
	\inp{u,u}_s=\norm{u}_s^2=\sum_{k\in\mathbb{Z}} |\inpo{u,\psi_k}|^2\theta_k^{2s},\qquad\forall\, u\in X_s.
\end{equation*}
The collection $\{X_s\}_{s\in\mathbb{Z}}$ is called Hilbert scale. 
\end{definition}
%
Note that $X_s$ is compactly embedded in $X_r$ if $s>r$, and is dense there.
Moreover, the spaces $X_s$ and $X_{-s}$
are conjugated with respect to the scalar product $\inpo{\cdot,\cdot}$,
that is
$$\norm{u}_s = \sup\{\inpo{u,v}:\, v\in X_{-s}\cap X_0, \norm{v}_{-s} = 1\}\qquad\forall\, u\in X_s\cap X_0,$$
so that the scalar product $\inpo{\cdot,\cdot}$ extends to a bilinear pairing $X_s\times X_{-s}\rightarrow \mathbb{R}$.
With a small abuse of notation we will denote this pairing as $X_0$-scalar product.

For the characterization of a symplectic structure on Hilbert scales we need to introduce the concept of morphism of order $d$.
\begin{definition}[Morphism of order $d$]
Let $s_0,s_1\in\mathbb{Z}$, $s_0\leq s_1$, be fixed and let $\{X_s\}_{s\in\mathbb{Z}}$ be a Hilbert scale.
Given a linear map $L:\bigcap_{s\in\mathbb{Z}} X_s\rightarrow \bigcup_{s\in\mathbb{Z}} X_s$ we denote by $\norm{L}_{\mathcal{L}(X_{s_1},X_{s_2})}$ its operator norm as a map $L:X_{s_1}\rightarrow X_{s_2}$. The map $L$ defines a morphism of order $d$ for $s\in[s_0,s_1]$ if $\norm{L}_{\mathcal{L}(X_{s},X_{s-d})}<\infty$ for any $s\in[s_0,s_1]$.
If the inverse map $L^{-1}$ exists and defines a morphism of order $-d$ for $s\in[s_0+d,s_1+d]$, then $L$ is an automorphism of order $d$ for $s\in[s_0,s_1]$.
The adjoint $L^*:X_{d-s}\rightarrow X_{-s}$ of a morphism of order $d$ for $s\in[s_0,s_1]$ is a morphism of order $d$ for $s\in[d-s_1,d-s_0]$.
\end{definition}
In particular, a bounded linear operator $L:X_{s_0}\rightarrow X_{s_0-d}$ can be regarded as a morphism of order $d$ for $s_0$.


Let $d\geq 0$ and $O_d$ be a non-empty open set in a Hilbert space $X_d$ from a Hilbert scale $\{X_s\}_{s\in\mathbb{Z}}$. For any $u\in O_d$ we identify the tangent space $T_uO_d$ with $X_d$, and treat differential $k$-forms as continuous functions
$$\alpha_u: \bigwedge{\!\!}^k\, X_d\longrightarrow \mathbb{R}$$
that are multi-linear and skew-symmetric.
In particular, 1-forms can be written as $a(u)\exd u(\xi):=\inpo{a(u),\xi}$, for any $\xi\in X_d$, where $a:O_d\rightarrow X_{-d}$ and $\exd$ is the exterior derivative; while
2-forms can be written as $\omega_u=A(u)\exd u\wedge\exd u$ such that
$\omega_u(\xi,\eta):=\inpo{A(u)\xi,\eta}$, for any $\xi,\eta\in X_d$, where
$A(u): X_d\rightarrow X_{-d}$ is a bounded self-adjoint operator.
This implies that $\omega_u:X_{\sigma}\times X_{\sigma}\rightarrow \mathbb{R}$ defines a continuous skew-symmetric bilinear form if $\sigma\geq - d/2$.

Given a $C^k$-smooth function $\Fcal:O_d\subset X_d\rightarrow \mathbb{R}$, $d\geq 1$, we consider the gradient map with respect to the pairing $\inp{\cdot,\cdot}$ as the $C^{k-1}$-smooth function
\begin{equation}
	\delta \Fcal:O_d\rightarrow X_{-d}\quad\mbox{such that}\quad\inp{\delta \Fcal(u),\xi}=\exd \Fcal(u)\xi\qquad\forall\,\xi\in X_d.
\end{equation}

A symplectic structure on the Hilbert scale can be introduced via the map
$\Jc(u): X_{s+d_J}\rightarrow X_s$, with $u\in O_d$,
a skew-adjoint morphism of order $d_J$ for $s\in[-d-d_J,d]$
called \emph{operator of the Poisson structure}.
\begin{definition}[Poisson bracket]\label{def:Pbra}
Let $\Fcal,\Lcal$ be $C^1$-smooth functions on $O_d$. Assume that they define gradient maps of orders $d_1$ and $d_2\leq 2d$ such that $d_1+d_2+d_J\leq 2d$. Then, the Poisson bracket of $\Fcal$ and $\Lcal$ is the continuous (on $O_d$) function
\begin{equation*}
	\bra{\Fcal}{\Lcal}(u)=\inpo{\Jc(u)\delta \Fcal(u),\delta \Lcal(u)}\qquad\forall\, u\in O_d.
\end{equation*}
\end{definition}
The Poisson bracket is skew-symmetric since $\Jc$ is skew-adjoint.
Moreover, it satisfies the Leibniz rule
$\bra{\Fcal\Lcal}{\mathcal{H}}=\bra{\Fcal}{\mathcal{H}}\Lcal+\Fcal\bra{\Lcal}{\mathcal{H}}$
and the Jacobi identity
$\bra{\Fcal}{\bra{\Lcal}{\mathcal{H}}}+\bra{\mathcal{H}}{\bra{\Fcal}{\Lcal}}+\bra{\Lcal}{\bra{\mathcal{H}}{\Fcal}}= 0$, for any $\Fcal,\Lcal,\mathcal{H}\in C^{2}(O_d)$.

The set $O_d$ endowed with the bracket $\bra{\cdot}{\cdot}$ forms a Poisson manifold, while the space $ C^{\infty}(O_d)$ of real-valued smooth functions over $(O_d,\bra{\cdot}{\cdot})$
together with the bracket $\bra{\cdot}{\cdot}$ forms a Lie algebra \cite[Proposition 3.3.17]{AM87}, called
the \emph{Poisson algebra} of $O_d$.

Hilbert scales endowed with a Poisson structure provide the scaffolding for the definition of Hamiltonian dynamics in infinite dimension.
Such dynamics
is described by a Hamiltonian vector field as follows.
By the bilinearity of $\bra{\cdot}{\cdot}$ and the Leibniz rule,
given a $C^1$-smooth function $\Hcal$ on $O_d$, 
there exists a locally unique vector field $V_{\Hcal}$ defined as the continuous map $V_{\Hcal}:O_d\rightarrow X_{-d-d_J}$
such that
$$\forall\, u\in O_d\qquad \Lie{V_{\Hcal}(u)}\Fcal=\bra{\Hcal}{\Fcal}(u),\qquad\forall\,\Fcal\in C^1(O_d),$$
where $\Lie{X}$ denotes the Lie derivative with respect to the velocity field $X$.
The vector $V_{\Hcal}(u)$ is called the \emph{Hamiltonian vector field} of $\Hcal$.
Let $\contr_{X}$ denote the contraction operator with the vector field $X$. Since, for any $\Fcal\in C^1(O_d)$,
$$\Lie{V_{\Hcal}(u)}\Fcal=\contr_{V_{\Hcal}(u)} \exd \Fcal=\exd \Fcal\,V_{\Hcal}(u)=\inpo{V_{\Hcal}(u),\delta \Fcal(u)}$$
and $\bra{\Hcal}{\Fcal}(u)=\inpo{\Jc(u)\dH(u),\delta \Fcal(u)}$, it follows that $V_{\Hcal}(u)=\Jc(u)\dH(u)$.

If the operator of the Poisson structure $\Jc(u): X_{s+d_J}\rightarrow X_s$ is an isomorphism, the operator $\overline{\Jc}(u):=-\Jc^{-1}(u):X_s\rightarrow X_{s+d_J}$ defines an skew-adjoint automorphism of order $-d_J$. Let us consider the closed 2-form $\omega_u:=\overline{\Jc}(u)\exd u\wedge\exd u$ where
$\overline{\Jc}(u):X_d\rightarrow X_{-d}$ depends $C^1$-smoothly on $u\in O_d$ and defines a linear isomoprhism $\overline{\Jc}(u):X_d\rightarrow X_{d+d_J}$ for any $d_J\geq 0$. The 2-form $\omega$ supplies $O_d$ with a symplectic structure and $\overline{\Jc}(u)$ is called operator of the symplectic structure.
The Hilbert scale $\{X_s\}_{s\in\mathbb{Z}}$ endowed with the 2-form $\omega$ forms a so-called \emph{symplectic Hilbert scale}, denoted as $(\{X_s\}_s,\omega)$.
The Poisson bracket on the symplectic space $(O_d,\omega)$ is defined as 
\begin{equation*}
	\bra{\Fcal}{\Lcal}(u)=\omega_u(V_{\Fcal},V_{\Lcal})\qquad\forall\, u\in O_d.
\end{equation*}

To a $C^1$-smooth function $\Hcal$ on $O_d$ the symplectic structure $\omega$ corresponds the Hamiltonian vector field $V_{\Hcal}$ defined as the continuous map
$V_{\Hcal}:O_d\rightarrow X_{-d-d_J}$ such that
\begin{equation*}
	\forall\, u\in O_d \qquad \omega_u(V_{\Hcal}(u),\xi)=-\exd \Hcal(u)\xi\qquad \forall\xi\in T_uO_d \simeq X_d.
\end{equation*}

With this definition, the Hamiltonian vector field $V_{\Hcal}$ satisfies
$\inpo{\overline{\Jc}(u) V_{\Hcal}(u),\xi}=-\inpo{\dH(u),\xi}$.
In particular, for any $C^1$-smooth function $\Hcal:O_d\times \mathbb{R}\rightarrow\mathbb{R}$, we denote by $V_{\Hcal}$ the non-autonomous vector field $V_{\Hcal}(u,t)=\Jc(u)\dH(u,t)$.
The corresponding Hamiltonian equation is
\begin{equation}\label{eq:HamEq}
\partial_t{u}(t)=V_{\Hcal}(u,t)=\Jc(u)\dH(u,t).
\end{equation}

A partial differential equation, supplemented by appropriate boundary conditions, is called a Hamiltonian PDE if, under a suitable choice of a symplectic Hilbert scale $(\{X_s\}_{s}, \omega)$, a domain $O_d\subset X_d$ and a Hamiltonian $\Hcal$, it can be written in the form \eqref{eq:HamEq}.

A vector field $V_{\Hcal}$ on a manifold $O_d$ determines a phase
flow, namely
a one-parameter group of diffeomorphisms
$\Phi^t_{V_{\Hcal}}:O_d\rightarrow O_d$ satisfying
$d_t\Phi^t_{V_{\Hcal}}(u)=V_{\Hcal}(\Phi^t_{V_{\Hcal}}(u))$
for all $t\in\mathcal{T}$ and $u\in O_d$, with $\Phi^0_{V_{\Hcal}}(u)=u$.

Hamiltonian dynamics is characterized by
the existence of differential invariants, and symmetry-related conservation laws.
\begin{definition}[Invariants of motion]\label{def:inv}
A function $\mathcal{I}\in C^{\infty}(O_d)$ is an \emph{invariant of motion} of 
the dynamical system with flow map
$\Phi^t_{V_{\Hcal}}$
if $\bra{\mathcal{I}}{\mathcal{H}}(u)=0$ for all $u\in O_d$.
Consequently, $\mathcal{I}$ is constant along the orbits of $V_{\Hcal}$.
\end{definition}

The Hamiltonian function, if time-independent, is an invariant of motion.
A particular subset of the invariants of motion of a dynamical system is
given by the \emph{Casimir invariants}, i.e. functions on $O_d$ which
$\bra{\cdot}{\cdot}$-commute with every other function on $O_d$.
\begin{definition}[Casimir invariants]\label{def:Cas}
If $\mathfrak{g}$ is a Lie algebra with Lie product $\bra{\cdot}{\cdot}{}$,
the \emph{centralizer} of a subset $\Ucal$ of $\mathfrak{g}$ is defined as
$\mathsf{C}_{\mathfrak{g}}(\Ucal):=\{\mathcal{C}\in\mathfrak{g}:\,\bra{\mathcal{C}}{\Fcal}{}=0\;\mbox{for all}\;\Fcal\in \Ucal\}$.
The centralizer $\mathsf{C}_{\mathfrak{g}}(\mathfrak{g})$ of the Lie algebra itself is called the \emph{centre} of
$\mathfrak{g}$ and its elements are called Casimir functions.
\end{definition}
%

\subsection{Gradient flow dissipation}\label{sec:Diss}
In this work we consider dissipation mechanisms in form of gradient flows. Since the seminal work \cite{jordan1998}, numerous contributions in the literature investigated the geometry and the behaviour of dissipative partial differential equations. Far from being exhaustive, we refer to \cite{santambrogio2017,mielke2011,arnrich2012} for a detailed discussion of the topic.






Leveraging the existence theory for gradient flows in Banach spaces, the Rothe method \cite{lions1969} in particular, we formulate the following gradient flow.
Let us consider a Hilbert scale $\{X_s\}_{s\in\mathbb{Z}}$. Note that each $X_s$ is reflexive and, since the injection $X_0\hookrightarrow X_s$, $s\geq 0$, is continuous, we can identify the Gelfand triple $X_s\subset X_0 \subset X_{-s}$.
%
%
Let us formulate a discrete in time gradient flow (whose limit is the continuous time gradient flow). Let $\mathcal{E}:X_s\rightarrow \mathbb{R}$ be a Fr\'{e}chet differentiable convex functional. Let $\tau>0$, and $v\in X_s$, we define
$$
u = \mathrm{arg} \inf_{\tilde{u}\in X_s} \Big(\frac{1}{2}\inpo{\tilde{u}-v,\tilde{u}-v} + \tau \mathcal{E}(\tilde{u})\Big).
$$
The functional augmented by the norm of the Hilbert space is clearly convex and Fr\'{e}chet differentiable. The Euler-Lagrange equations read as follows:
$$
\inpo{u-v + \tau \delta\mathcal{E}(u), z} = 0\qquad \forall z \in X_s. 
$$
Since the above expression holds for any $z\in X_s$, we take $z= u-v$ and get
$$
\| u-v\|_0^2 + \tau \inpo{\delta\mathcal{E}(u), u-v} = 0,
$$
which gives
\begin{equation}
\label{eq:nablaE}
\inpo{\delta\mathcal{E}(u), v-u} = \frac{1}{\tau} \| u-v\|_0^2 \geq 0.
\end{equation}
Let us now exploit the convexity of the functional $\mathcal{E}$ to show that the energy is dissipated. 
The convexity of $\mathcal{E}$ gives, for $\vartheta \in [0,1]$ and $u,v\in X_s$, that
$\mathcal{E}\left( (1-\vartheta)u + \vartheta v \right) \leq (1-\vartheta)\mathcal{E}(u) + \vartheta \mathcal{E}(v)$.
This entails
$
\vartheta^{-1}\left[ \mathcal{E}\left( u + \vartheta (v-u) \right) - \mathcal{E}(u)\right] \leq \mathcal{E}(v) - \mathcal{E}(u)
$ for any $\vartheta \in [0,1]$.
In the limit $\vartheta\rightarrow 0$,
$$
\inpo{\delta \mathcal{E}(u), v-u} \leq \mathcal{E}(v) - \mathcal{E}(u).
$$
Using \eqref{eq:nablaE} yields
$$
0 \leq \frac{1}{\tau} \norm{u-v}_{0}^2 \leq \mathcal{E}(v) - \mathcal{E}(u),
$$
which clearly shows that $\mathcal{E}(u)\leq\mathcal{E}(v)$.
The continuous gradient flow is the limit, for $\tau\rightarrow 0$, of the time discrete gradient flow and it enjoys the same property of energy dissipation. Its weak formulation can be stated as follows: for all $v\in X_s$
$$
\inpo{\partial_t u, v} + \inpo{\delta\mathcal{E}(u), v} = 0.
$$
It is possible to show \cite{lions1969} that there exists a unique solution $u\in L^p(\mathcal{T},V)$, $1\leq p\leq \infty$, of this problem in the temporal interval $\Tcal$.

In the present work we consider gradient flows with respect to a  (potentially singular) metric. Let $\mathcal{S}=-\mathcal{E}$ be a concave functional and $\Gc(u)$ be a positive semi-definite self-adjoint operator. This gives rise to the following formulation:
$$
\partial_t u = \Gc(u)\dS(u).
$$
Let us consider the case in which $\Gc(u)$ has an inverse $\Gc(u)^{-1}$. In such a case, the gradient flow can be considered as the continuous limit of the following discrete in time gradient flow
$$
u = \mathrm{arg} \inf_{\tilde{u}\in X_s} \Big(\frac{1}{2}\inpo{\tilde{u}-v, \Gc(v)^{-1} (\tilde{u}-v)} - \tau \mathcal{S}(\tilde{u})\Big).
$$
We can show that the quantity $\mathcal{S}$ is non-decreasing.

Let us consider the case in which $\mathcal{S}: X_0 \rightarrow \mathbb{R}$ is the entropy, and $\Gc(u): X_{-s} \rightarrow X_{s}$ is a linear self-adjoint operator. The equation for the time evolution of the entropy reads
$$
d_t \mathcal{S}(u) = \inpo{\dS(u), \Gc(u) \dS(u)}\geq 0.
$$
It is possible to introduce a positive semi-definite bilinear form $(\cdot,\cdot):X_{s}\times X_{s}\rightarrow\mathbb{R}$ defined as
\begin{equation}\label{eq:Gprod}
\big(\Fcal,\mathcal{L}\big)(u)=\inpo{\Gc(u)\delta \Fcal(u),\delta \mathcal{L}(u)}\qquad\forall\, u\in X_s,\, \Fcal,\mathcal{L}:X_s\rightarrow\mathbb{R}.
\end{equation}

We present some examples of dissipation terms, obtained as a gradient flow. 

\begin{example}[Heat equation]
Let $\Omega\subset\mathbb{R}^n$ be a Lipschitz continuous domain.
We consider the heat flow in the following setting:
$X_0=L^2(\Omega)$ and $u\in H^1_0(\Omega) = X_1$.
By Rellich-Kondrakov theorem, the injection $X_1\hookrightarrow X_0$ is compact.
Let $\Gc(u)=-\Delta$ with domain $H^1_0(\Omega)$ and let
$$
\Scal(u) = -\frac{1}{2}\int_{\Omega} u^2 \ dx.
$$
Since its Fr\'{e}chet derivative is $\dS(u) = -u$ one gets
$$
d_t \mathcal{S}(u) = -\int_{\Omega} u \Delta u \ dx = \int_{\Omega} \nabla u \cdot \nabla u \ dx.
$$
We obtain, as entropy evolution, the classical evolution of the $L^2$ norm of the solution of the heat equation. 

The classical formulation of the heat equation as a gradient flow is the following. Let $X_0 = L^2(\Omega)$ and $X_1=H^1_0(\Omega)$ as above. The entropy is the negative Dirichlet energy:
$$
\mathcal{S}(u) = -\frac{1}{2}\int_{\Omega} |\nabla u|^2 \ dx.
$$
Its Fr\'{e}chet derivative is $\delta\mathcal{S}(u) = \Delta u \in X_{-1}=X_1'$. The operator $\Gc(u)$ is the identity, and for all $z\in H^1_0(\Omega)$:
$$
\inpo{\Delta u, z} = - \inpo{\nabla u, \nabla z}.
$$
Remark that when we take $z=u$ we get: $\inpo{\dS(u), u} \leq 0$.
\end{example}

\begin{example}[Double bracket dissipation]
Another example of dissipation structure is provided by the so-called double bracket, see e.g. \cite{bloch1996,sato2024}. In this formulation, the dissipation operator is given by $\Gc = -\Jc^2$ or $\Gc = -\overline{\Jc}^2$, where $\Jc$ and $\overline{\Jc}$ are the Poisson and symplectic operators defined in \Cref{sec:Ham}. 
Let us show that the operator $\Gc= -\Jc^2$ is positive definite and self-adjoint. Let $v\in X_{s+d_J}$, it holds:
$$
\inpo{\Gc(u) v, v} = -\inpo{\Jc^2(u) v, v}
= -\inpo{\Jc(u) v, \Jc^*(u) v}
= \inpo{\Jc(u) v, \Jc(u) v} \geq 0.
$$
The self-adjointness follows by observing that, for any $v,\xi\in X_{s+d_J}$,
$$
-\inpo{\Jc^2(u) v, \xi}
= -\inpo{\Jc(u) v, \Jc^*(u) \xi}
= \inpo{\Jc(u) v, \Jc(u) \xi}
= -\inpo{\Jc^*(u) v, \Jc(u) \xi}
= -\inpo{v, \Jc^2(u) \xi}.
$$
This is enough to show that the operator $\Jc^2$ induces an increase in the entropy functional and, therefore, it represents a candidate to describe dissipation mechanisms. Remark that, if we can consider, as the entropy, $\Scal=-\Hcal$, the double bracket dissipates the total energy (Hamiltonian) of the system.
\end{example}

\subsection{Hamiltonian systems with gradient flow dissipation}
The focus of this work is on Hamiltonian partial differential equations that are supplied by a gradient flow dissipation.
The setting is as follows:
let us consider the Gelfand triple $X_{s}\subset X_0\subset X_{-s}$
and let $\Tcal:=(0,T]\subset\mathbb{R}$ be a bounded temporal interval.
We consider the evolution equation \eqref{eq:pbm}, namely: find $u\in L^p(\Tcal,X_s)$ such that
\begin{equation*}
\partial_t u = \Jc(u)\dH(u) + \Gc(u)\dS(u),\qquad t\in \Tcal,
\end{equation*}
supplied by the initial condition $u(0)=u_0\in X_s$
and suitable boundary conditions that we assume to be encoded in the definition of the space $X_s$.
Moreover, we assume that the functions $\mathcal{H},\mathcal{S}: X_s\rightarrow \mathbb{R}$ 
satisfy $\mathcal{H}(u), \mathcal{S}(u)< +\infty$ for any $u\in X_s$. 
Further, the Poisson structure $\Jc$ and the operator $\Gc$ satisfy the conditions of \Cref{sec:Ham,sec:Diss}
and are such that $\Jc(u)\dH\in X_{-s}$ and $\Gc(u)\dS \in X_{-s}$ for any $u\in X_s$.

With a small abuse of notation, we will refer to the function $\Scal$ appearing in the generic evolution equation \eqref{eq:pbm} as ``entropy'' even if the physical interpretation of such quantity might differ based on the particular example considered. 

We present some examples of Hamiltonian systems with a dissipative term.

\begin{example}[Advection-diffusion equation]
Let us consider the advection-diffusion equation
\begin{equation*}
\partial_t u + \mathbf{v}\cdot \nabla u - \Delta u=0
\end{equation*}
in a Lipschitz continuous spatial domain $\Omega\subset\mathbb{R}^n$, with constant velocity $\mathbf{v}$ and homogeneous Dirichlet boundary conditions.
Let $X_0=L^2(\Omega)$ and $X_1=H^1_0(\Omega)$. 
Let us define $\mathcal{H}:X_1\rightarrow \mathbb{R}$ as
$$
\mathcal{H}(u) = \frac{1}{2}\int_{\Omega} |u|^2 \ dx.
$$
Then, the Fr\'{e}chet derivative of $\mathcal{H}$ is $\delta\mathcal{H}(u) = u$.
We also introduce the operators $\Gc(u)=\Gc=\Delta$, with domain $X_1$ and $\Jc(u)=\Jc=-\mathbf{v}\cdot\nabla$. 
It holds
$\Jc\delta\mathcal{H}(u)=-\mathbf{v}\cdot\nabla u\in X_0\subset X_{-1}$,
and
$\Gc\delta\mathcal{H}(u)=\Delta u\in X_{-1}$.
Thus, the advection-diffusion problem can then be cast as
\begin{equation*}
\partial_t u=(\Jc+\Gc)\delta\Hcal(u).
\end{equation*}
The kinetic energy of the system dissipates. Indeed,
$$
d_t \Hcal(u) = \int_{\Omega} u \Delta u \ dx = -\int_{\Omega} |\nabla u|^2 \ dx.
$$
\end{example}

\begin{example}[Isentropic KdV equation]
Let us consider the Korteweg-de Vries equation
\begin{equation}\label{eq:kdv}
    \partial_t u + \alpha u\partial_x u+\eta \partial_x^3 u = \nu\partial_x^2 u
\end{equation}
in a one-dimensional spatial domain $\Omega\subset\mathbb{R}$ with periodic boundary conditions. We assume that $\alpha,\eta,\nu$ are positive bounded constants.

Let  $X_0=L^2(\Omega)$ and $X_1=H^1(\Omega)$, we introduce the quantities
\begin{equation*}
\Jc(u)=\Jc=-\partial_x,\qquad
\Hcal(u) = \int_{\Omega}\bigg(\dfrac{\alpha}{6}u^3-\dfrac{\eta}{2}(\partial_x u)^2\bigg) dx.
\end{equation*}
The Fr\'{e}chet derivative of the Hamiltonian reads
$$
\dH(u) = \frac{\alpha}{2} u^2 + \eta \partial^2_x u.
$$
For the KdV equation in Gardner's bracket formulation \cite{Gar71}, a possible entropy is given by
\begin{equation}\label{eq:KdVmass}
    \Scal(u) = \int_{\Omega} u \ dx.
\end{equation}
This can be interpreted as the mass of the wave (especially when considering solitons). When having periodic boundary conditions, or a closed system, it is therefore natural to consider an evolution equation which preserves the mass, $d_t \Scal = 0$, and dissipates an energy.

Let the operator $\Gc(u):X_{1}\rightarrow X_{-1}$ be defined as
$\Gc(u)v=\nu v\partial_x^2 u$ for any $v\in X_{1}$.
Then, the KdV equation \eqref{eq:kdv} can be cast as in \eqref{eq:pbm}.
In this setting the entropy is conserved since it holds
$$
d_t \Scal(u) = \inpo{\dS(u), \partial_t u} =
\inpo{1, \Jc \dH(u)} + \nu \inpo{1, \partial^2_x u} = 0,
$$
in view of the periodic boundary conditions and the fact that
$1\in \mathrm{ker}(\Jc)$.

The kinetic energy of the system
\begin{equation*}
    \mathcal{E}(u)=\int_{\Omega} \dfrac{u^2}{2}\, dx
\end{equation*}
instead dissipates since
$$
d_t \Ecal(u) = \inpo{\delta \Ecal(u), \partial_t u} =
\inpo{u, -\alpha u\partial_x u-\eta\partial_x^3 u} + \nu \inpo{u, \partial^2_x u} = - \nu \norm{\partial_x u}_{0}^2\leq 0.
$$
Although the Hamiltonian evolution does not have a definite sign, it is possible to show that what dissipates is the Hamiltonian of the so-called boosted soliton, that is the energy given by the Hamiltonian augmented by the kinetic energy of a reference frame moving at a certain velocity $\beta>0$, that is
\begin{equation*}
    \widehat{\Ecal}(u):=\Hcal(u)+\beta\int_{\Omega} \dfrac{u^2}{2}\, dx.
\end{equation*}
It is possible to show that there exists $\beta_0$ such that, for all $\beta>\beta_0$, we have $d_t \widehat{\Ecal}<0$. Indeed,
$$
d_t \widehat{\Ecal}(u) = \inpo{\delta \widehat{\Ecal}(u), \partial_t u}
= \inpo{\dH(u), \partial_t u} + \beta \inpo{u,\partial_t u}
= \nu\eta \norm{\partial^2_x u}_{0}^2 +\dfrac{\nu\alpha}{2} \inpo{\partial_x^2 u,u^2}-\nu\beta\norm{\partial_x u}_{0}^2. 
$$

Moreover, an equilibrium solution is given by
$$
u_* = \varrho_0=\frac{1}{|\Omega|}\int_{\Omega} u \ dx,
$$
and the mass \eqref{eq:KdVmass} is simply $\Scal =\varrho_0 |\Omega|$.
Indeed, we have
$$
\dH|_{u_*} = \frac{\alpha}{2}\varrho_0^2,\quad\mbox{and}\quad \partial^2_x u_0 = 0,
$$
and, hence, 
$$
\partial_t u_* = \dfrac{\alpha}{2}\Jc\varrho_0^2 = 0, 
$$
since the constant function is in the kernel of the operator $\Jc$. 

\end{example}

\subsubsection{The case of metriplectic or GENERIC systems}
\label{sec:metr}

A special subclass of problems that can be described via \eqref{eq:pbm}
is the one where the dynamics is generated by both a Hamiltonian and an entropy function. The corresponding geometric
formalism has been introduced under the name of metriplectic structure \cite{PJM86} in plasma physics and
GENERIC formalism \cite{GEN97} in the context of non-equilibrium thermodynamics.

Metriplectic systems are characterized by an evolution equation as in \eqref{eq:pbm},
and further equipped with so-called \emph{mutual degeneracy conditions}, namely
\begin{equation}\label{eq:HScond}
    \Gc(u)\delta\Hcal(u) = 0
    \quad\mbox{and}\quad
    \Jc(u)\delta\Scal(u) = 0
    \qquad\forall\, u\in X_s.
\end{equation}
Note that, in view of the skew-adjointness of the operator $\Jc$ and the self-adjointness of $\Gc$, these conditions imply that $\bra{\Hcal}{\Scal}(u)=(\Hcal,\Scal)(u)=0$ for any $u\in X_s$. This, in turn,
ensures that the Hamiltonian $\Hcal$ is a conserved quantity of motion while the entropy $\Scal$ increases in time. Indeed,
\begin{equation*}
\begin{aligned}
& d_t\Hcal(u) =
        \inpo{\delta\Hcal(u),\partial_t u}
        =\bra{\Hcal}{\Hcal}(u) + (\Hcal,\Scal)(u) = 0.\\
& d_t\Scal(u) =
        \inpo{\delta\Scal(u),\partial_t u}
         = \bra{\Hcal}{\Scal}(u) + (\Scal,\Scal)(u)
       \geq 0.
\end{aligned}
\end{equation*}

\section{The variational formulation of the problem}
\label{sec:weakform}

The variational formulation of the evolution problem \eqref{eq:pbm} reads: Given $u_0\in X_s$, find $u\in L^p(\mathcal{T},X_s)$ such that
$$
\inpo{\partial_t u, v} = \inpo{\Jc(u)\dH(u), v} + \inpo{\Gc(u)\dS(u), v},\qquad \forall\,v\in X_s.
$$
Since $\Jc(u)$ is skew-adjoint and $\Gc(u)$ is self-adjoint, we can rewrite the weak formulation as
\begin{equation}
\label{eq:weak_formulation}
\inpo{\partial_t u, v} = -\inpo{\dH(u), \Jc(u) v} + \inpo{\dS(u), \Gc(u) v},\qquad\forall\, v\in X_s.
\end{equation}
Equation \eqref{eq:weak_formulation} is the starting point to derive a mixed variational formulation, detailed in the next section.

\subsection{Mixed formulation}\label{sec:mixed}
Let $s,d,e\in \mathbb{N}\setminus\{0\}$. Let us consider the Hilbert spaces $X_{s-d}$ and $X_{s-e}$ and remark that $X_s$ is dense in $X_{s-d}$ and in $X_{s-e}$. In the following, we are going to assume that $\delta\mathcal{H}\in X_{-s+d}$ and $\delta\mathcal{S}\in X_{-s+e}$. 

The mixed formulation reads: Given $u_0\in X_s$, find $(u,z,y)\in C^1(\mathcal{T};X_s)\times C^0(\mathcal{T};X_{-s+d}) \times C^0(\mathcal{T}; X_{-s+e})$ such that
\begin{equation}\label{eq:mixed}
\begin{cases}
\inpo{\partial_t u, v} = -\inpo{z, \Jc(u) v} + \inpo{y, \Gc(u)v},& \qquad \forall\, v\in X_s, \\
\inpo{z-\delta\mathcal{H}, v}= 0,& \qquad\forall\, v\in X_s, \\
\inpo{y-\delta\mathcal{S}, v}= 0,& \qquad\forall\, v\in X_s.
\end{cases}
\end{equation}

\subsubsection{Equivalence of the variational and mixed formulations}

Let us consider the subspaces $K_{\Jc}, K_{\Gc}\subset X_s$ defined as
\begin{equation*}
\begin{aligned}
& K_{\Jc}:=\{u\in X_s:\; \Jc(u)v\in X_{s-d} \;\mbox{for all}\,v\in X_s\},\\
& K_{\Gc}:=\{u\in X_s:\; \Gc(u)v\in X_{s-e} \;\mbox{for all}\,v\in X_s\}.
\end{aligned}
\end{equation*}


\begin{lemma} 
\label{lem:density}
Let $\xi\in X_{-s+d}$ and $u\in K_{\Jc}$. If $\inpo{\xi, v} = 0$ for all $v\in X_s$, then
$$
\inpo{\xi, \Jc(u) v} = 0,\qquad\forall\, v\in X_s.
$$
\begin{proof}
Let $\tilde{v}\in X_s$. Adding and subtracting this element to the second component of the inner product gives
$$
|\inpo{\xi,\Jc(u) v}|=|\inpo{\xi, \tilde{v}}
+\inpo{\xi, \Jc(u) v - \tilde{v}}|.
$$
The first term of the right-hand side vanishes by hypothesis, since $\tilde{v}\in X_s$. The second term can be bounded, using H\"{o}lder's inequality, as
$$
|\inpo{\xi,\Jc(u)v-\tilde{v}}| \leq \|\xi\|_{-s+d}\|\Jc(u)v-\tilde{v}\|_{s-d}.
$$
Since $X_s$ is dense in $X_{s-d}$ we can choose $\tilde{v}$ in such a way that the second term 
becomes arbitrarily small, and the conclusion follows.
\end{proof}
\end{lemma}
Remark that, by using the same density argument, we can prove that, given $\zeta\in X_{-s+e}$ and $u\in K_{\Gc}$, if $\inpo{\zeta, v} = 0$ for all $v\in X_s$ then $\inpo{\zeta, \Gc(u) v}= 0$.

\begin{proposition}
Assume that $\dH\in X_{-s+d}$ and $\dS\in X_{-s+e}$.
If $t\in\Tcal\mapsto u(t)\in X_s$ solution of the mixed formulation \eqref{eq:mixed} belongs to $K_{\Jc}\cap K_{\Gc}\subset X_s$, then the mixed formulation is equivalent to the weak formulation \eqref{eq:weak_formulation} of the evolution equation \eqref{eq:pbm}.
\begin{proof}
By virtue of the result of \Cref{lem:density},
the solution $(u(t),z(t),y(t))\in X_s\times X_{-s+d} \times X_{-s+e}$ of the mixed formulation \eqref{eq:mixed} satisfies, for any $t\in\Tcal$, 
$\inpo{z-\dH, \Jc(u) v} = 0$, and $\inpo{y-\dS, \Gc(u) v} = 0$ for any $v\in X_s$.
From the first equation of \eqref{eq:mixed}, we thus get, for any $v\in X_s$,
$$\inpo{\partial_t u, v} = -\inpo{z, \Jc(u) v}+\inpo{y, \Gc(u) v}=-\inpo{\dH(u), \Jc(u) v}+\inpo{\dS(u), \Gc(u) v}.$$
\end{proof}
\end{proposition}

Remark that, under the assumption that $X_s$ is dense in both $X_{-s+d}$ and $X_{-s+e}$, that is $2s\leq d,e$, then $\inpo{\xi,v}=0$ for any $v\in X_s$ directly implies $\xi=0\in X_{-s+d}$. A similar result holds for $X_{-s+e}$.

\section{The semi-discrete problem}
\label{sec:semidisc}

\subsection{Conformal variational discretisation}\label{sec:confdisc}
Under the assumption that $X_s$ is dense in both $X_{-s+d}$ and $X_{-s+e}$ $(2s \geq d, e)$, we build a conformal variational discretisation of the mixed formulation \eqref{eq:mixed} by introducing a finite dimensional subspace $V_{N}\subset X_s$, closed in $X_s$, and such that, for $N\rightarrow\infty$, $V_N$ is dense in $X_s$. Let $N\in\mathbb{N}^*$ and $\left\lbrace v_i \right\rbrace_{i=1}^N$ be a basis of $V_{N}$. We restrict the mixed formulation by taking test functions which belong to $V_N$. Moreover, by considering a Galerkin method, we look for numerical approximations $u_N,z_N,y_N$ in $V_N$. Let $a,b,c:\mathcal{T}\rightarrow \mathbb{R}^N$, we write the numerical approximation of the solution, for any $(t,x)\in \mathbb{R}^+\times \Omega$, as:
$$
u_N(t,x) = \sum_{i=1}^N a_i(t) v_i(x),\qquad
z_N(t,x) = \sum_{i=1}^N b_i(t) v_i(x),\qquad
y_N(t,x) = \sum_{i=1}^N c_i(t) v_i(x).
$$
This leads to the following discretisation of the mixed formulation \eqref{eq:mixed}:
\begin{equation}\label{eq:mixed_disc}
\left\{
\begin{aligned}
&\sum_{j=1}^N \langle v_j, v_i \rangle_0\, d_t a_j = 
- \sum_{j=1}^N \big( \langle v_j, \Jc(u_N) v_i \rangle_0 b_j + \langle v_j, \Gc(u_N) v_i \rangle_0 c_j\big), \ \ \ & 1\leq i\leq N, \\
&\sum_{j=1}^N \langle v_j, v_i \rangle_0 b_j - \langle \delta\mathcal{H}(u_N), v_i \rangle_0 = 0, \ \ \ & 1\leq i\leq N, \\
& \sum_{j=1}^N \langle v_j, v_i \rangle_0 c_j - \langle \delta\mathcal{S}(u_N), v_i \rangle_0 = 0, \ \ \ & 1\leq i\leq N.
\end{aligned}\right.
\end{equation}

In order to rewrite the discrete formulation of the problem, let us state the following result.
\begin{lemma}
\label{lem:discrete_Ham_Ent}
Let $u_N\in V_N$ and let $a\in\mathbb{R}^N$ be the coefficients representing $u_N$ in the basis $\left\lbrace v_i \right\rbrace_{i=1}^N$ of $V_N$. Let $H,S:\mathbb{R}^N\rightarrow\mathbb{R}$ be defined as
$$
H(a) = \mathcal{H}(u_N) = \mathcal{H}\left( \sum_{i=1}^N a_i v_i \right),\quad\mbox{and}
\quad
S(a) = \mathcal{S}(u_N) = \mathcal{S}\left( \sum_{i=1}^N a_i v_i \right).
$$
Then, the Fr\'{e}chet derivatives of $\mathcal{H}$ and $\mathcal{S}$, evaluated in $u_N \in V_N$, are such that
$$
\partial_{a_i} H(a) = \langle \delta\mathcal{H}(u_N), v_i \rangle_0,
\qquad
\partial_{a_i} S(a) = \langle \delta\mathcal{S}(u_N), v_i \rangle_0,\qquad\forall\,1\leq i\leq N.
$$
\begin{proof}
The proof is performed by direct computation starting from the definition of the Fr\'{e}chet derivative. For all $v\in X_s$, it holds
$$
\langle \delta\mathcal{H}(u), v \rangle_0 = \lim_{\vartheta\rightarrow 0} \frac{\mathcal{H}(u+\vartheta v)- \mathcal{H}(u)}{\vartheta}.
$$
When $u=u_N\in V_N$, the variation is restricted to $v\in V_N,\ v= \sum_{j=1}^N q_j v_j$, and we can write
$$
\langle \delta\mathcal{H}(u_N), v \rangle_0 = \lim_{\vartheta\rightarrow 0} \frac{\mathcal{H}(u_N+\vartheta v)- \mathcal{H}(u_N)}{\vartheta}= \lim_{\vartheta\rightarrow 0} \frac{\mathcal{H}\left(\sum_{j=1}^N (a_j + \vartheta q_j)v_j \right)- \mathcal{H}(u_N)}{\vartheta}.
$$
Since $\left\lbrace q_j \right\rbrace_{1\leq j\leq N}$ is arbitrary, we can take $q_j = \delta_{ij}$ and get:
$$
\langle \delta\mathcal{H}(u_N), v_i \rangle_0 = \lim_{\vartheta\rightarrow 0} \frac{\mathcal{H}\left(\sum_{j=1}^N (a_j + \vartheta \delta_{ij})v_j \right)- \mathcal{H}\left(\sum_{j=1}^N a_j v_j\right )}{\vartheta} = \partial_{a_i} H(a).
$$
The proof for the Fr\'{e}chet derivative of $\mathcal{S}$ is analogous.
\end{proof}
\end{lemma}
Let $u_N\in V_N$ and let $a\in\mathbb{R}^N$ be the coefficients of the expansion of $u_N$ in the basis $\left\lbrace v_i \right\rbrace_{i=1}^N$ of $V_N$. 
We introduce the following matrices: the mass matrix $M\in\mathbb{R}^{N\times N}$, and $\Jd(a)\in\mathbb{R}^{N\times N}$ and $\Gd(a)\in \mathbb{R}^{N\times N}$ which are the representation of the symplectic and metric operators, respectively. Their components read, for any $1\leq i,j \leq N$,
\begin{equation}\label{eq:MJG}
M_{i,j} = \inpo{v_j, v_i},\qquad 
\Jd(a)_{i,j} = \inpo{\Jc(u_N) v_j, v_i},\qquad
\Gd(a)_{i,j} = \inpo{\Gc(u_N) v_j, v_i}.
\end{equation}
With this notation and the results of \Cref{lem:discrete_Ham_Ent}, the semi-discrete (in space) mixed formulation of the problem can be written, for any $t\in \Tcal$, as:
\begin{equation}\label{eq:mixed_disc_1}
\begin{cases}
M d_t a(t) = \Jd(a(t)) b(t) + \Gd(a(t)) c(t)\\
M b(t) - \nabla H(a(t)) = 0, \\
M c(t) - \nabla S(a(t)) = 0.
\end{cases}
\end{equation}

\subsection{Geometric structure of the semi-discrete problem}
In the next Lemma we prove that the discretisation of the continuous Poisson bracket from \Cref{def:Pbra} yields a bilinear form that is still a Poisson bracket.
\begin{lemma}
\label{lem:poisson_bracket}
Let 
\begin{equation}\label{eq:Jhat}
\Jm(a) := M^{-1}\Jd(a)M^{-1}\qquad\forall\, a \in\mathbb{R}^N.
\end{equation}
Let the bilinear form $\bra{\cdot}{\cdot}_N:C^1(\mathbb{R}^N) \times C^1(\mathbb{R}^N) \rightarrow C^1(\mathbb{R}^N)$ be defined as:
$$
\bra{F}{L}_N(a) = \Jm(a) \nabla F(a)\cdot  \nabla L(a),\qquad\forall\, a \in\mathbb{R}^N,
$$
where $\cdot$ denotes the Euclidean product in $\mathbb{R}^N$ and
$F,L:\mathbb{R}^N\rightarrow \mathbb{R}$ are real valued functions in $C^1(\mathbb{R}^N)$.
Let $\Fcal,\Lcal: X_s \rightarrow \mathbb{R}$ be $C^1$-functionals, such that, for all elements $u_N\in X_s$ with $u_N = \sum_{i=1}^N a_i v_i$, it holds $\Fcal(u_N) = F(a)$ and $\Lcal(u_N) = L(a)$. Then,
\begin{enumerate}
\item $\bra{F}{L}_N$ is the restriction to $V_N\subset X_s$ of the continuous Poisson bracket $\bra{\Fcal}{\Lcal}$.
\item $\bra{\cdot}{\cdot}_N$ is a Poisson bracket, namely it is skew-symmetric, it satisfies the Leibniz rule and the Jacobi identity.
\end{enumerate}
\begin{proof}
Let $\Fcal, \Lcal:X_s\rightarrow\mathbb{R}$
and $u=u_N\in V_N$. Then
$\delta \Fcal(u_N)=\sum_{j=1}^N f_j v_j$ where
$f_j$, by virtue of \Cref{lem:discrete_Ham_Ent}, is the $j$th component of the vector $f=M^{-1}\nabla F(a)$. Similar considerations can be made for $\delta \Lcal$. Using \Cref{def:Pbra} of the Poisson bracket, it holds
\begin{equation*}
    \begin{aligned}
        \{\Fcal,\Lcal\}(u_N)
        & =\inpo{\Jc(u_N)\delta \Fcal(u_N),\delta \Lcal(u_N)}
=\sum_{i,j=1}^N\inpo{\Jc(u_N)f_jv_j,\ell_i v_i}\\
& =\sum_{i,j,k,m=1}^N M^{-1}_{ik}M^{-1}_{jm}\inpo{\Jc(u_N)v_j,v_i} \partial_{a_k}L(a)\partial_{a_m}F(a) \\
& = \sum_{k,m=1}^N \big(\nabla L\big)_k(a) \left( \sum_{i,j=1}^N M^{-1}_{k,i} \Jd(a)_{i,j} M^{-1}_{j,m} \right) \big(\nabla F\big)_m(a) =\bra{F}{L}_N(a).
    \end{aligned}
\end{equation*}
It is easy to show that $\bra{\cdot}{\cdot}_N$ is skew-symmetric and satisfies the Leibniz rule.
The second item follows from the fact that
the continuous bracket satisfies the Jacobi identity.
\end{proof}
\end{lemma}
Observe that an analogous result holds for the bilinear form associated with $\Gc$.
Indeed, the bilinear form \eqref{eq:Gprod} restricted to $V_N$ gives
\begin{equation*}
        \big(\Fcal,\mathcal{L}\big)(u_N)
         =\inpo{\Gc(u_N)\delta \Fcal(u_N),\delta \mathcal{L}(u_N)}
= \sum_{k,m=1}^N \big(\nabla L\big)_k(a) \left( \sum_{i,j=1}^N M^{-1}_{k,i} \Gd(a)_{i,j} M^{-1}_{j,m} \right) \big(\nabla F\big)_m(a). 
\end{equation*}
The discrete positive semi-definite bilinear form associated with \eqref{eq:Gprod} is then
 \begin{equation*}
     \big(F,L\big)_N(a) = \Gm(a) \nabla F(a)\cdot\nabla L(a),\qquad
     \Gm(a):=M^{-1}\Gd(a) M^{-1},\quad\forall\, a\in\mathbb{R}^N.
 \end{equation*}

We are now in the position of proving the main results of the present work. In particular, we are going to show that the semi-discrete in space mixed formulation \eqref{eq:mixed_disc_1} is a finite dimensional Hamiltonian system with a dissipation mechanism in the form of a finite dimensional gradient flow. We are going to separate the proof in two propositions, treating the cases $\Jc(u)=0$ and $\Gc(u) = 0$ separately.
\begin{proposition}\label{prop:Ham}
Let $\Gc(u) = 0$. Then, the mixed formulation \eqref{eq:mixed_disc_1}:
\begin{enumerate}
\item It is a consistent semi-discretised form of the continuous Hamiltonian system obtained by setting $\Gc(u)=0$ in \eqref{eq:pbm}. 
\item It is Hamiltonian, with the discrete Poisson bracket defined in \Cref{lem:poisson_bracket}.
\end{enumerate}
\begin{proof}
Since $\Gc(u)=0$ we can restrict to the first and second equations of the mixed formulation \eqref{eq:mixed_disc_1} and take $c=0$. 
From the second equation we have
$b = M^{-1}\nabla H(a)$, that we can substitute
into the first equation and get
$$
d_t a  = M^{-1}\Jd(a) M^{-1} \nabla H(a)=\Jm(a)\nabla H(a),
$$
and this concludes the first part of the proof. 

Moreover, \Cref{lem:poisson_bracket} ensures that
$\Jm(a)$ is a discrete Poisson tensor.
%
This is enough to conclude that the system is a finite dimensional Hamiltonian system.
\end{proof}
\end{proposition}

Let us now prove that the discretisation of the gradient flow leads to a finite dimensional gradient flow. 
\begin{proposition}\label{prop:diss}
Let $\Jc(u) = 0$. Then, the mixed formulation \eqref{eq:mixed_disc_1}:
\begin{enumerate}
\item It is a consistent semi-discretised form of the continuous gradient flow system obtained by setting $\Jc(u)=0$ in \eqref{eq:pbm}. 
\item It is a finite dimensional gradient flow for the entropy $S$.
\end{enumerate}
\begin{proof}
Since $\Jc(u)=0$ we can simply set $b=0$ and consider the first and third equations in the mixed formulation \eqref{eq:mixed_disc_1}.
Replacing
$c = M^{-1} \nabla S(a)$
into the mixed formulation yields
$$
d_t a = M^{-1}\Gd(a) M^{-1} \nabla S(a)=\Gm(a)\nabla S(a).
$$
Let us remark that $\Gd(a)$ is a symmetric positive semi-definite matrix by construction and $M$ is symmetric and positive definite. Then, the discrete dissipation tensor $\Gm(a)$ is symmetric positive semi-definite. This implies that the discrete system entropy $S(a)$ cannot decrease:
$$
d_t S(a) = \nabla S^{\top}(a) \Gm(a) \nabla S(a) \geq 0,
$$
and this concludes the proof.
\end{proof}
\end{proposition}

\subsection{Conservation properties of the semi-discrete problem}

Another interesting question concerns the quantities which are conserved at discrete level by the mixed discretisation we have described. Concerning the Hamiltonian part, a major role is played by the Casimirs of the system, as introduced in \Cref{def:Cas}.

\begin{proposition}
Let $u\in X_s$ and $\Gc(u) = 0$.
The following holds:
\begin{enumerate}
\item If, for any $u\in X_s$, $\mathrm{ker}(\Jc(u))$ has finite dimension $r\leq N$, and $\mathrm{ker}(\Jc(u)) \subset V_N$, then, for any $a\in\mathbb{R}^N$, $\Jm(a)$ has a kernel of dimension $r$, and the discretisations of the continuous Casimirs are discrete Casimirs. 
\item If, for any $u\in X_s$, $\mathrm{ker}(\Jc(u)) \perp_{X_s} V_N$, then the restriction of any Casimir to $V_N$ is a constant function. 
\end{enumerate}
\begin{proof}
Let $\Fcal:X_s\rightarrow\mathbb{R}$ be a Casimir of $\Jc$ that is
$$\bra{\Fcal}{\Ccal}(u)=\inpo{\Jc(u)\delta\Fcal(u),\delta\Ccal(u)}=0\qquad\forall\, u\in X_s,\;\Ccal:X_s\rightarrow\mathbb{R}.$$
If $\mathrm{ker}(\Jc(u)) \subset V_N$ for any $u\in X_s$, then the Fr\'{e}chet derivative $\delta\Fcal$ of the Casimir $\Fcal$ belongs to $V'_N$. This means that we can write 
$$\delta \Fcal(u_N) = \sum_{j=1}^N f_j v_j\qquad\forall\, u_N=\sum_{j=1}^N a_j v_j,$$
where $f_j=f_j(a) = (M^{-1}\nabla F(a))_j$, and $F(a) = \Fcal(u_N)$.

Moreover, since $\Fcal$ is a Casimir, it holds
$$
0=\inpo{\Jc(u_N) \delta \Fcal(u_N),v} = \sum_{j=1}^N f_j(a) \inpo{\Jc(u_N) v_j,v},\qquad\forall\, u_N\in V_N,\, v\in X_s.$$
Taking $v=v_i$, we get $(\Jd(a)f(a))_i=0$.
This implies that the vector $f(a)=M^{-1}\nabla F(a)$ belongs to the kernel of $J(a)$, for any $a\in\mathbb{R}^N$, and, hence, $\nabla F(a)$ belongs to the kernel of $\Jm(a)=M^{-1}\Jd(a)M^{-1}$.
This concludes the first part of the proof.

Concerning the second item, let us proceed by contradiction. Let us assume that the discrete Casimir $F$ depends on $a$ so that $\nabla F(a)$ is, in general, different from zero. It follows that there exists a function $\delta\Fcal\in V_N$ such that
its coefficients in the basis $\{v_i\}_{i=1}^N$ are given by
$f_i = \sum_{l=1}^N M^{-1}_{il} \partial_{a_l} F(a)$. This function is the Fr\'{e}chet derivative of a functional $\Fcal$ restricted to $V_N$. For $\Fcal$ to be a Casimir, it must hold: $\delta\Fcal \in \mathrm{ker}(\mathcal{J}(u))$. This could be possible only if at least part of the kernel is in $V_N$, hence contradicting the hypothesis. 
\end{proof}
\end{proposition}



\begin{proposition}[Conservation of invariants]
\label{prop:invariants}
Let $\Kcal\subseteq X_s$ be the phase space of problem \eqref{eq:pbm}. Assume that $\Kcal$ is convex and contains the origin.
Let $\Fcal: \Kcal\rightarrow \mathbb{R}$ be an invariant of motion of problem \eqref{eq:pbm} as in \Cref{def:inv}.
Let $F:\mathbb{R}^N\rightarrow \mathbb{R}$ be a real valued function belonging to $C^1(\mathbb{R}^N)$ defined as $F(a)=\Fcal(u_N)$ for any $u_N=\sum_{i=1}^N a_i v_i\in V_N\subset \Kcal$.
Then, $F$ is an invariant of motion of the semi-discrete problem \eqref{eq:mixed_disc_1}.
\begin{proof}
Let $u_N\in V_N$. Since $\Fcal$ is an invariant of \eqref{eq:mixed} it holds
$\bra{\Fcal}{\Hcal}(u_N)+(\Fcal,\Scal)(u_N)=0$.
This is equivalent to
$\bra{F}{H}_N(a)+(F,S)_N(a)=0$. Since $u_N$ is chosen arbitrarily, $F$ is an invariant of motion of the
semi-discrete problem \eqref{eq:mixed_disc_1}.
\end{proof}
\end{proposition}

Next result shows that the equilibria of \eqref{eq:pbm}
belonging to the approximation space $V_N$ are equilibria of the semi-discrete problem.
\begin{proposition}
Let $u^*\in X_s$ be an equilibrium point of problem \eqref{eq:pbm} that is
$$\Jc(u^*)\dH(u^*) + \Gc(u^*)\dS(u^*)=0.$$
If $u^*$ belongs to $V_N$, then it
is equilibrium point of the semi-discrete system \eqref{eq:mixed_disc}.
    \begin{proof}
Since $u^*\in V_N$ it can expressed as
$u^*=\sum_{i=1}^N a_i^* v_i$.
Moreover, since $u^*$ is equilibrium point it holds
$\inpo{f(u^*),v}=0$ for any $v\in X_s$, where
$f(u):=\Jc(u)\dH(u) + \Gc(u)\dS(u)$ for any $u\in X_s$.
Taking $v=v_i\in V_N$ yields
\begin{equation*}
    \begin{aligned}
      0=\inpo{f(u^*),v_i}
      & =\sum_{j=1}^N b_j\inpo{\Jc(u^*)v_j,v_i}
      +\sum_{j=1}^N c_j\inpo{\Gc(u^*)v_j,v_i}\\
      & =\sum_{j=1}^N b_j\Jd(a^*)_{i,j}
      +\sum_{j=1}^N c_j\Gd(a^*)_{i,j}
      = \big(\Jd(a^*)b+\Gd(a^*)c\big)_i.
    \end{aligned}
\end{equation*}
Since $b=M^{-1}\nabla H(a^*)$ and $c=M^{-1}\nabla S(a^*)$, one has
$\Jd(a^*)M^{-1}\nabla H(a^*)+\Gd(a^*)M^{-1}\nabla S(a^*)=0$ and the conclusion follows.
    \end{proof}
\end{proposition}

\subsubsection{The case of metriplectic or GENERIC systems}

The mixed variational formulation \eqref{eq:mixed}
holds also for the metriplectic systems introduced in \Cref{sec:metr}. The mutual degeneracy conditions are imposed weakly as
\begin{equation}\label{eq:HScondweak}
    \inpo{\Gc(u)\delta\Hcal(u),w} = 0
    \quad\mbox{and}\quad
    \inpo{\Jc(u)\delta\Scal(u),w} = 0
    \qquad\forall\, u,w\in X_s.
\end{equation}
The semi-discretisation of problem \eqref{eq:mixed}
in $V_N\subset X_s$ can be performed as in \eqref{eq:mixed_disc_1} for metriplectic systems as well.
Moreover, the mutual degeneracy conditions \eqref{eq:HScondweak} are discretised by imposing $w=v_i\in V_N$ and taking $u_N\in V_N$; this yields
\begin{equation}\label{eq:semidiscr_nullcond}
    0 = \inpo{\Gc(u_N)\delta\Hcal(u_N),v_i}
    = \sum_{j,k=1}^N M^{-1}_{jk}\partial_{a_k} H(a) \inpo{\Gc(u_N) v_j,v_i}
    = (\Gd(a)M^{-1}\nabla H(a))_i,
\end{equation}
and analogously $\Jd(a)M^{-1}\nabla S(a)=0$
for all $a\in \mathbb{R}^N$.
The resulting semi-discrete system remains metriplectic, as shown in the next result.
\begin{proposition}
The semi-discrete in space mixed formulation \eqref{eq:mixed_disc_1} together with the conditions \eqref{eq:semidiscr_nullcond} is a finite
dimensional metriplectic system.
\begin{proof}
Following the reasoning in the proofs of \Cref{prop:Ham,prop:diss}, the semi-discrete system reads
$$d_t a(t) = \Jm(a(t)) \nabla H(a(t)) + \Gm(a(t)) \nabla S(a(t)).$$
Moreover, the discretisation of the mutual degeneracy conditions \eqref{eq:semidiscr_nullcond}  yields
\begin{equation}\label{eq:discr_nullcond}
    \Gm(a)\nabla S(a)=0\qquad\mbox{and}\qquad\Jm(a)\nabla S(a)=0.
\end{equation}
%
%
\end{proof}
\end{proposition}
As a consequence of the previous result, the Hamiltonian $H$ is a conserved quantity and $S$ is increasing. Indeed, from \eqref{eq:discr_nullcond}, it holds
\begin{equation*}
    \begin{aligned}
        & (H,S)_N(a)=\nabla H(a)^{\top}\Gm(a)\nabla S(a) =0,\\
        & \bra{H}{S}_N(a)=\nabla H(a)^{\top}\Jm(a)\nabla S(a)=0.
    \end{aligned}
\end{equation*}
The results follows from the skew-symmetry of $\Jm(a)$ and the positive semi-definiteness of $\Gm(a)$, for any $a\in\mathbb{R}^N$.


\subsection{Convergence results}
Considering the convergence properties of the proposed semi-discretisation, we first notice that, in all the cases in which the mixed formulation is equivalent to the classical variational formulation, it enjoys the same convergence properties. 
In this section, we propose a proof of convergence in norm $\| \cdot \|_0$ which makes use of Lipschitz continuity hypotheses of the Fr\'{e}chet derivatives and of the operators $\Jc$ and $\Gc$. In particular, let us consider the following assumptions:
\begin{enumerate}
\item [$(H_1)$] The Fr\'{e}chet derivatives $\delta \mathcal{H}$ and $\delta\mathcal{S}$ are Lipschitz continuous: there exist $c_H,c_S>0$ such that
$$
\| \delta\mathcal{H}(u)-\delta\mathcal{H}(v)\|_{-s+d} \leq c_H \| u-v\|_{s},\quad\mbox{and}
\quad \| \delta\mathcal{S}(u)-\delta\mathcal{S}(v)\|_{-s+e} \leq c_S \| u-v\|_{s}.
$$
for any $u,v\in X_s$.
\item [$(H_2)$] The map $u\in X_s\mapsto \Jc(u)\in\mathcal{L}(X_{-s+d},X_{-s})$
is Lipschitz continuous: there exists a constant $c_J>0$ such that, for any $u,v\in X_s$,
$$
\| \Jc(u) - \Jc(v) \|_{\mathcal{L}(X_{-s+d},X_{-s})} \leq c_J \| u-v\|_{s}.
$$
Analogously, concerning the dissipation, we assume that the map 
$u\in X_s\mapsto \Gc(u)\in\mathcal{L}(X_{-s+e},X_{-s})$
is Lipschitz continuous: there exists a constant $c_G>0$ such that, for any $u,v\in X_s$,
$$
\| \Gc(u) - \Gc(v) \|_{\mathcal{L}(X_{-s+e},X_{-s})} \leq c_G \| u-v\|_{s}.
$$
\end{enumerate}

\begin{proposition}[Convergence of the semi-discretisation scheme]
Let $s,d,e\in \mathbb{N}\setminus\{0\}$, $2s\geq d,e$. 
Under the hypotheses $(H_1)$ and $(H_2)$, the solution $u_N(t)\in V_N$ of the semi-discrete system  \eqref{eq:mixed_disc_1} converges to the solution of the continuous mixed formulation \eqref{eq:mixed}: For any  $t\in\mathcal{T}$,
$$
\lim_{N\rightarrow\infty}\norm{u(t)-u_N(t)}_{0} = 0.
$$

\begin{proof}
In what follows, for the sake of compactness of the notation, we omit the time dependence on $u(t),z(t),y(t)$ and write simply $u,z,y$.

Let $P_N$, $Q_N$, and $R_N$ be three residual functions, orthogonal to $V_N$, that is, $\langle P_N, v_N \rangle_0 = 0$, $\langle Q_N, v_N \rangle_0 = 0$, and $\langle R_N, v_N \rangle_0 = 0$, for any $v_N\in V_N$, obtained by testing the numerical approximation of the mixed problem \eqref{eq:mixed_disc} against a generic element of the space, $v\in X_s$:
\begin{equation}\label{eq:mixed_disc_Vtest}
  \begin{cases}
\inpo{\partial_t u_N, v} = \inpo{\Jc(u_N)z_N + \Gc(u_N)y_N + Q_N, v},\\
\inpo{z_N - \dH(u_N)-P_N, v} = 0,\\
\inpo{y_N - \dS(u_N)- R_N, v} = 0.
\end{cases}
\end{equation}
Note that, since, for $N\to\infty$, $V_N$ is dense in $X_s$, it holds $\norm{Q_N}_{-s}\to 0$ as $N\to\infty$ and similarly for $R_N$ and $P_N$ in the associated norms.
Let the error be the function $e = u-u_N$ where $u$ is solution of the mixed formulation \eqref{eq:mixed} and $u_N$ is solution of \eqref{eq:mixed_disc_Vtest}. We derive the mixed formulation for the error by taking the difference of the mixed formulation for the solution $u$ and the one for its numerical approximation $u_N$:
for any $v\in X_s$,
\begin{equation*}
\begin{cases}
\inpo{\partial_t e, v} = \inpo{\Jc(u) z + \Gc(u)y -\Jc(u_N)z_N - \Gc(u_N)y_N -  Q_N, v},\\
\inpo{z-z_N - \dH(u) + \dH(u_N) + P_N, v} = 0,\\
\inpo{y - y_N - \dS(u) + \dS(u_N) + R_N, v}= 0.
\end{cases}
\end{equation*}
As a particular test function, we choose $v=e$ and get
\begin{equation*}
\begin{cases}
\inpo{\partial_t e, e} = \inpo{\Jc(u) z + \Gc(u)y -\Jc(u_N)z_N - \Gc(u_N)y_N - Q_N, e},\\
\inpo{z-z_N - \dH(u) + \dH(u_N) + P_N, e} = 0,\\
\inpo{y - y_N - \dS(u) + \dS(u_N) + R_N, e}= 0.
\end{cases}
\end{equation*}
From the first equation, by making use of H\"{o}lder's inequality, we can write:
$$
\frac{1}{2}\partial_t \norm{e}_{0}^2 \leq
\norm{\Jc(u) z + \Gc(u)y -\Jc(u_N)z_N - \Gc(u_N)y_N}_{-s}\norm{e}_{s} + \norm{Q_N}_{-s}\norm{e}_{s}.
$$
We focus on the first term on the right-hand side.
By means of the triangular inequality we can separate the Hamiltonian and dissipative parts, 
$$
\norm{\Jc(u) z + \Gc(u)y -\Jc(u_N)z_N - \Gc(u_N)y_N}_{-s} \leq
\norm{\Jc(u) z - \Jc(u_N)z_N}_{-s} + \norm{\Gc(u)y - \Gc(u_N)y_N}_{-s}.
$$
We can treat the terms separately. Concerning the first term, we add and subtract $\Jc(u)z_N$ and make use of the triangular inequality:
$$
\norm{\Jc(u) z - \Jc(u_N)z_N}_{-s} \leq
\norm{\Jc(u) (z - z_N)}_{-s} + \norm{\Jc(u)-\Jc(u_N)) z_N}_{-s}.
$$
Using the Lipschitz continuity of $\Jc(u)$ with respect to $u$ (hypothesis $(H_2)$), we write:
\begin{equation}\label{eq:tmp_deltaJ}
\norm{\Jc(u) z - \Jc(u_N)z_N}_{-s} \leq
\widehat{c}_J(u) \norm{z - z_N}_{-s+d} + c_J \norm{e}_{s}\norm{z_N}_{-s+d},
\end{equation}
where $\widehat{c}_J(u):=\norm{\Jc(u)}_{\mathcal{L}(X_{-s+d},X_{-s})}$.
In order to bound $\norm{z-z_N}_{-s+d}$ we consider the equation for the element of the cotangent bundle
$$
\inpo{z-z_N, v} = \inpo{\delta \mathcal{H}(u) - \delta \mathcal{H}(u_N) - P_N, v}\qquad\forall\, v\in X_s.
$$
Let us consider the norm in $X_{-s+d}$ as given by
$$
\norm{z}_{-s+d} = \sup_{z_*\in X_{s-d}\setminus\{0\}} \frac{|\inpo{z, z_*}|}{\norm{z_*}_{s-d}}.
$$
Let $v_*\in X_{s-d}$, $\norm{v_*}_{s-d}=1$, be supremiser of the duality product $\inpo{z-z_N, v}$. 
Since $X_s\subset X_{s-d}$ with a dense inclusion,
we consider the sequence in $X_s$ that converges to $v^*$.
We get:
$$
\norm{z-z_N}_{-s+d} = 
|\inpo{\dH(u) - \dH(u_N) - P_N, v_*}| \leq
\norm{\dH(u) - \dH(u_N)}_{-s+d}
+ \norm{P_N}_{-s+d}.
$$
The Hamiltonian is Fr\'{e}chet differentiable, and, by virtue of the hypothesis $(H_1)$ we have Lipschitz continuity of the Fr\'{e}chet derivative. This makes it possible to introduce a constant $c_H>0$ such that:
$$
\norm{z-z_N}_{-s+d} \leq c_H \norm{e}_{s}  + \norm{P_N}_{-s+d}.
$$
We use this result into \eqref{eq:tmp_deltaJ} and get:
$$
\norm{\Jc(u) z - \Jc(u_N)z_N}_{-s} \leq
\widehat{c}_J(u) \big(c_H\norm{e}_{s} + \norm{P_N}_{-s+d}\big) + c_J \norm{e}_{s}\norm{z_N}_{-s+d},
$$
which can be rewritten as
$$
\norm{\Jc(u) z - \Jc(u_N)z_N}_{-s} \leq \big(c_H \widehat{c}_J(u)+ c_J\norm{z_N}_{-s+d}\big)\norm{e}_{s} +
\widehat{c}_J(u)\norm{P_N}_{-s+d}.
$$
Concerning the dissipation mechanism, we can proceed in the same way, and get
$$
\norm{\Gc(u) y - \Gc(u_N)y_N}_{-s} \leq
\big(c_S \widehat{c}_G(u)
 + c_G\| y_N \|_{-s+e}\big)\norm{e}_{s} +
 \widehat{c}_G(u)\norm{R_N}_{-s+e},
$$
where $\widehat{c}_G(u):=\| \Gc(u) \|_{\mathcal{L}(X_{-s+e},X_{-s})}$.
Summing these two terms, we get that there exists $C_N(u)>0$ such that
\begin{equation*}
    \norm{\Jc(u) z - \Jc(u_N)z_N + \Gc(u) y - \Gc(u_N)y_N}_{-s} \leq C_N(u) \norm{e}_{s}
     + \widehat{c}_J(u)\norm{P_N}_{-s+d}+
    \widehat{c}_G(u)\norm{R_N}_{-s+e},
\end{equation*}
where $C_N(u):=c_H \widehat{c}_J(u)+ c_J\norm{z_N}_{-s+d} + c_S \widehat{c}_G(u) + c_G\| y_N \|_{-s+e}$.
Henceforth, it holds
$$
\frac{1}{2}\partial_t \norm{e}_{0}^2 \leq
C_N(u) \norm{e}_{s}^2 + \big(\norm{Q_N}_{-s} + \| \Jc(u) \|_{\mathcal{L}(X_{-s+d},X_{-s})}\norm{P_N}_{-s+d} + \| \Gc(u) \|_{\mathcal{L}(X_{-s+d},X_{-s})}\norm{R_N}_{-s+e} \big)\norm{e}_{s}.
$$
%
The equation is homogeneous in $\| e\|_{s}$ and since we can approximate the initial condition of the system, $e_0=u(0)-u_N(0)$ is such that $\lim_{N\rightarrow \infty}\| e_0\|_{s} = 0$, we can conclude that the mixed semi-discrete formulation converges to the continuous mixed formulation of the system. 
Furthermore, the error with respect to the solution of the problem tends to zero in the norm of $X_0$. 
\end{proof}
\end{proposition}


\section{Temporal discretisation of the semi-discrete problem}
\label{sec:time}

Concerning the numerical time integration of the semi-discrete problem \eqref{eq:mixed_disc_1} we consider two well-established methods: the Average Vector Field (AVF) introduced in \cite{QMacL08} and the implicit midpoint rule.

First, we consider a partition of the temporal interval $\Tcal=(0,T]$ into $N_t$ subintervals $(t^k,t^{k+1}]$, $0\leq k\leq N_t-1$. For simplicity of exposition, we assume this partition to be uniform so that $t^{k}=t^0+k\dt$, for $0\leq k\leq N_t$ with $\dt= N_t^{-1}T$.
The extension to the non-uniform case is trivial.

Given the ordinary differential equation
\begin{equation}\label{eq:ode}
    d_t a = f(a),\qquad a(0)=a^0,
\end{equation}
the AVF scheme approximates the solution at time $t^{k+1}$, $k\geq 0$, as
\begin{equation}\label{eq:AVF}
     \dfrac{a^{k+1}-a^k}{\dt} = \int_0^1 f\big((1-\xi)a^k+\xi a^{k+1}\big)\,d\xi.
\end{equation}
In our case, the velocity field of the flow is given by
$$
f(a) = \Jm(a) \nabla H(a) + \Gm(a)\nabla S(a). 
$$
\begin{lemma}
Assume that $\Jd$ and $\Gd$ in \eqref{eq:mixed_disc_1} are constant-valued. Then, the AVF scheme applied to \eqref{eq:mixed_disc_1}
\begin{itemize}
    \item exactly preserves the Hamiltonian, when $\Gd=0$;
    \item satisfies the dissipative law, when $\Jd=0$;
    \item is a metriplectic time integrator.
\end{itemize}
\begin{proof}
First note that, if $\Jd$ and $\Gd$ are constant-valued then so are
$\Jm$ and $\Gm$.
Let $y(\xi):=(1-\xi)a^k+\xi a^{k+1}$ and let $\cdot$ denote the Euclidean inner product. Multiplying Equation \eqref{eq:AVF} by the integral of $\nabla H(y(\xi))$ in $[0,1]$, one gets
\begin{equation*}
\begin{aligned}
\dfrac{1}{\dt}\int_0^1 (a^{k+1}-a^k)\cdot\nabla H(y(\xi))\,d\xi=
&\,\int_0^1 \nabla H(y(\xi))\,d\xi\, \Jm
\int_0^1 \nabla H(y(\xi))\,d\xi\\
& +\int_0^1 \nabla H(y(\xi))\,d\xi\, \Gm
\int_0^1 \nabla S(y(\xi))\,d\xi.
\end{aligned}
\end{equation*}
If $\Gd=0$ or $\Gd\nabla H(y)=0$, for any $y\in \mathbb{R}^N$,
then the above equality gives
\begin{equation*}
    \dfrac{1}{\dt} \big(H(a^{k+1})-H(a^k)\big)=0
\end{equation*}
 owing to the skew-symmetry of $\Jm$.
A similar reasoning can be carried out for the dissipation of $S$, leading to the conclusion.
\end{proof}
\end{lemma}
Note that, the AVF scheme, in its classic version, is a second order temporal integrator. It is, however, not a symplectic integrator.

The implicit midpoint rule applied to the ODE \eqref{eq:ode} is given by
\begin{equation}\label{eq:IM}
     \dfrac{a^{k+1}-a^k}{\dt} = f\big(a^{k+1/2}\big),
\end{equation}
where the midpoint state $a^{k+1/2}$ is defined as $a^{k+1/2}=\frac{a^{k+1}+a^k}{2}$ for any $0\leq k\leq N_t$.

\begin{lemma}\label{lem:IMconsHam}
    Assume that $H$ is a quadratic function, i.e., $H(a)=a^{\top}La$ with $L$ symmetric positive definite.
    Then, the implicit midpoint rule applied to problem \eqref{eq:mixed_disc_1} ensures exact conservation of the Hamiltonian whenever $\Gd=0$ or the system is metriplectic.
\begin{proof}
By the linearity of $\nabla H$ one has, for any $k\geq 0$,
    \begin{equation*}
    \begin{aligned}
        \dfrac{1}{\dt} \big(H(a^{k+1}) - & H(a^k)\big) =  \dfrac{1}{\dt} \nabla H(a^{k+1/2})\cdot(a^{k+1}-a^k) = \nabla H(a^{k+1/2})\cdot f(a^{k+1/2})\\
        = &\, \nabla H(a^{k+1/2})\cdot \Jm(a^{k+1/2})\nabla H(a^{k+1/2})
         + \nabla H(a^{k+1/2})\cdot \Gm(a^{k+1/2})\nabla S(a^{k+1/2})\\
        = &\,\nabla H(a^{k+1/2})\cdot \Gm(a^{k+1/2})\nabla S(a^{k+1/2}).
\end{aligned}
\end{equation*}
The last term vanishes when either $\Gd=0$ or the null condition $\Gd\nabla H(y)=0$ holds for any $y\in \mathbb{R}^N$.
\end{proof}
\end{lemma}
An analogous reasoning holds for the dissipation of $S$.
\begin{lemma}
    Assume that $S$ is a quadratic function, i.e., $S(a)=a^{\top}La$ with $L$ symmetric positive definite.
    Then, the implicit midpoint rule applied to problem \eqref{eq:mixed_disc_1} gives dissipation of $S$
    whenever $\Jd=0$ or the system is metriplectic.
\end{lemma}
The result of the previous Lemmas implies that, if $H$ and $S$ are quadratic functions of the state, then the implicit midpoint rule is a metriplectic integrator.
Note that, for linear equations, the AVF scheme and the implicit midpoint rule coincides.

\section{Numerical tests}
\label{sec:num}
In this section we present several numerical tests to illustrate the properties of the method. 
In all numerical tests we take $X_0=L^2(\Omega)$ with the $L^2$ inner product denoted as $\inpo{\cdot,\cdot}$.

The errors will be measured in relative $L^2$ norm, and we define:
$$
e(t) = \frac{\| u(t) - u_N(t) \|_{L^2(\Omega)}}{\| u(t) \|_{L^2(\Omega)}}.
$$
We will take the $L^{\infty}(\mathcal{T})$ norm in time:
$$
\| e(t) \|_{L^{\infty}(\mathcal{T})} = \sup_{t\in\mathcal{T}} e(t),
$$
and, in the following, we will refer to it as $L^{\infty}(\mathcal{T},L^2(\Omega))$ relative error.

For the spatial approximation of the problem we consider a finite element partition $\Omega_N$ of the domain $\Omega$ in the sense of Ciarlet. In the numerical tests of the next sections $V_N$ will be the space of piecewise linear finite elements on the partition $\Omega_N$ and $N$ will correspond to the number of mesh nodes.

\subsection{Korteveg-de Vries (KdV) equation}
The KdV equation is an example of non-linear dispersive wave equation, which is an infinite dimensional integrable Hamiltonian system \cite{Gar71}.
Let $\Omega = [0,20\pi]$, we will consider here its standard form, with a dissipation in a form of a heat term, as presented in Eq. \eqref{eq:kdv},
$$
\partial_t u + \alpha u \partial_x u + \eta \partial^3_x u = \nu \partial^2_x u,
$$
with periodic boundary conditions. In what follows, we consider $\alpha = 6$, $\eta = 1$, and we ley $\nu$ vary. First we consider the non-dissipative case, which corresponds to the classical KdV equation, and then the dissipative one. Let us recall that the Hamiltonian is:
$$
\mathcal{H}(u) = \int_{\Omega} \frac{\alpha}{6} u(x)^3 - \frac{\eta}{2} (\partial_x u(x))^2 \ dx,
$$
and $\Jc(u) = \Jc = -\partial_x$, which does not depend on the solution. We consider the entropy
$$
\mathcal{S}(u) = -\frac{1}{2}\int_{\Omega} u(x)^2 \ dx, 
$$
with $\Gc(u) = \Gc = - \nu \partial^2_x$.

In order to choose the semi-discretisation in space, let us write the mixed formulation by assuming the functions regular enough and check the minimal regularity required for the formulation to hold. Let $v\in H^1_{\mathrm{per}}(\Omega)$ where the subscript refers to the periodic boundary conditions. We have
\begin{equation*}
    \begin{cases}
        \langle \partial_t u, v \rangle_0 = -\langle \partial_x z, v \rangle_0 + \nu \langle \partial_x y, \partial_x v \rangle_0, \\
        \langle z, v\rangle_0 = \langle \frac{\alpha}{2} u^2, v \rangle_0 - \eta \langle \partial_x u , \partial_x v \rangle_0,\\
        \langle y, v \rangle_0 = -\langle u, v \rangle_0.
    \end{cases}
\end{equation*}
Let us remark that, since the entropy is minus the energy of the solution $\delta \mathcal{S} = -u$, the variable $y$ is somehow redundant. However, we keep it here in order to study the mixed formulation of the problem. We observe that choosing $V_N$ such that $V_N\subset H^1_{\mathrm{per}}(\Omega)$ provides enough regularity to make the formulation wellposed. We use piecewise linear finite elements $V_N = \mathbb{P}^1_{\mathrm{per}}(\Omega_N)$ on a uniform partition $\Omega_N$ with $N$ nodes. 

The mass 
$$
\mathcal{M}(u) = \int_{\Omega} u(x) \ dx
$$ 
is a Casimir of the Hamiltonian part of the system; indeed $\delta \mathcal{M} = 1 \in \mathrm{ker}(\Jc)$. The dissipative term, given the periodic boundary conditions, does not affect the mass. We can therefore conclude that the mass is a conserved quantity of the system.

Given the results of \Cref{prop:invariants}, since $1\in \mathbb{P}^1_{\mathrm{per}}(\Omega_N)$, we can infer that the proposed semi-discretisation conserves the mass. 

Let $M,K\in\mathbb{R}^{N\times N}$ be the mass and stiffness matrices respectively, whose components are defined as:
$$M_{i,j}:=\inpo{v_j,v_i}\qquad K_{i,j}:=\inpo{\partial_x v_j,\partial_x v_i}.$$
Using the notation from \Cref{sec:confdisc}, the semi-discretisation of the system reads as follows
\begin{equation}\label{eq:KdVsemi-disc}
    \begin{cases}
        M d_t a=- Ab-\nu Ka\\
        M b = \dfrac{\alpha}{2} Ra\otimes a-\eta K a
    \end{cases}
\end{equation}
where $A_{i,j}:=\inpo{\partial_x v_j,v_i}$ and $R\in\mathbb{R}^{N\times N^2}$ is the triple tensor defined as
$$R_{k,h}=\inpo{v_jv_k,v_i}\qquad h=(i-1)N+j,\quad\forall\; 1\leq i,j\leq N.$$
For the numerical time integration of the semi-discrete problem \eqref{eq:KdVsemi-disc} we consider the implicit midpoint rule and the AVF method described in \Cref{sec:time}.
In both discretisations the linear terms are evaluated at the midpoint state, while the nonlinear terms are discretised differently according to the chosen integrator.
In particular, the fully discrete system corresponding to \eqref{eq:KdVsemi-disc} reads
\begin{equation*}
    \begin{cases}
     M \dfrac{a^{k+1}-a^k}{\dt} = -A b^* - \nu Ka^{k+1/2},\\[2ex]
    M b^* = \alpha R(d_2^{-1} a^{k+1}\otimes a^{k+1}+d_2^{-1} a^{k}\otimes a^{k}+d_1^{-1} a^{k+1}\otimes a^{k})-\eta Ka^{k+1/2}.
    \end{cases}
\end{equation*}
It can be easily verified that, for the implicit midpoint rule, $b^*=b^{k+1/2}$ which corresponds to $d_1=4$, $d_2=8$; whereas
for the AVF, $d_1=d_2=6$.
We consider, as unknown, the vectors $a^{k+1}, b^{*}$, which are the discretisation of the elements of the manifold and its cotangent bundle, respectively. Moving the unknowns to the left hand side, we get
\begin{equation*}
\left\{\begin{aligned}
& \left(M+\frac{\dt}{2}\nu K\right) a^{k+1} + \dt A b^* = \left(M-\frac{\dt}{2}\nu K\right) a^{k}\\
& \left(\frac{\eta}{2}K -  d_1^{-1} \alpha R a^{k}\otimes \right)a^{k+1} - d_2^{-1}\alpha R a^{k+1} \otimes a^{k+1} + M  b^{*} = \left(-\frac{\eta}{2} K + d_2^{-1}\alpha R a^{k} \otimes\right) a^{k}.
\end{aligned}\right.
\end{equation*}
The nonlinearity in the latter equation is treated by means of a fixed point iteration.

In the numerical tests of this subsection we set the number $N$ of degrees of freedom of the spatial discretisation to $N=512$ and consider $N_t=800$ equispaced time iterations in the temporal interval $\Tcal=[0,15]$.

\subsubsection*{Classical, non dissipative solutions}
In this section, we consider $\nu=0$, and get the classical KdV equation. We are in the periodic boundary conditions setting, but since the domain is large, the one soliton solution on the real line is close to the true solution of the system, and we will use it as a benchmark. 

The analytical solution is:
\begin{equation}\label{eq:soliton}
u(x,t) = \mathrm{sech}^2\Bigg( \frac{\sqrt{2}}{2} (x - 5\pi -2t) \Bigg),
\end{equation}
which corresponds to a solitary wave propagating at constant speed.

\begin{figure}[H]
\includegraphics[scale=0.5]{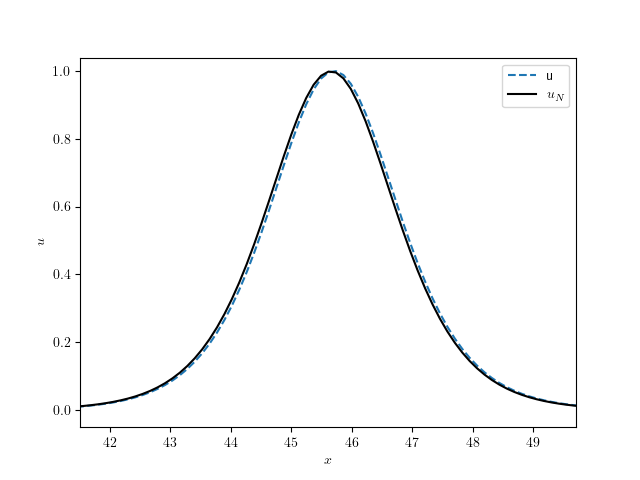}
\includegraphics[scale=0.5]{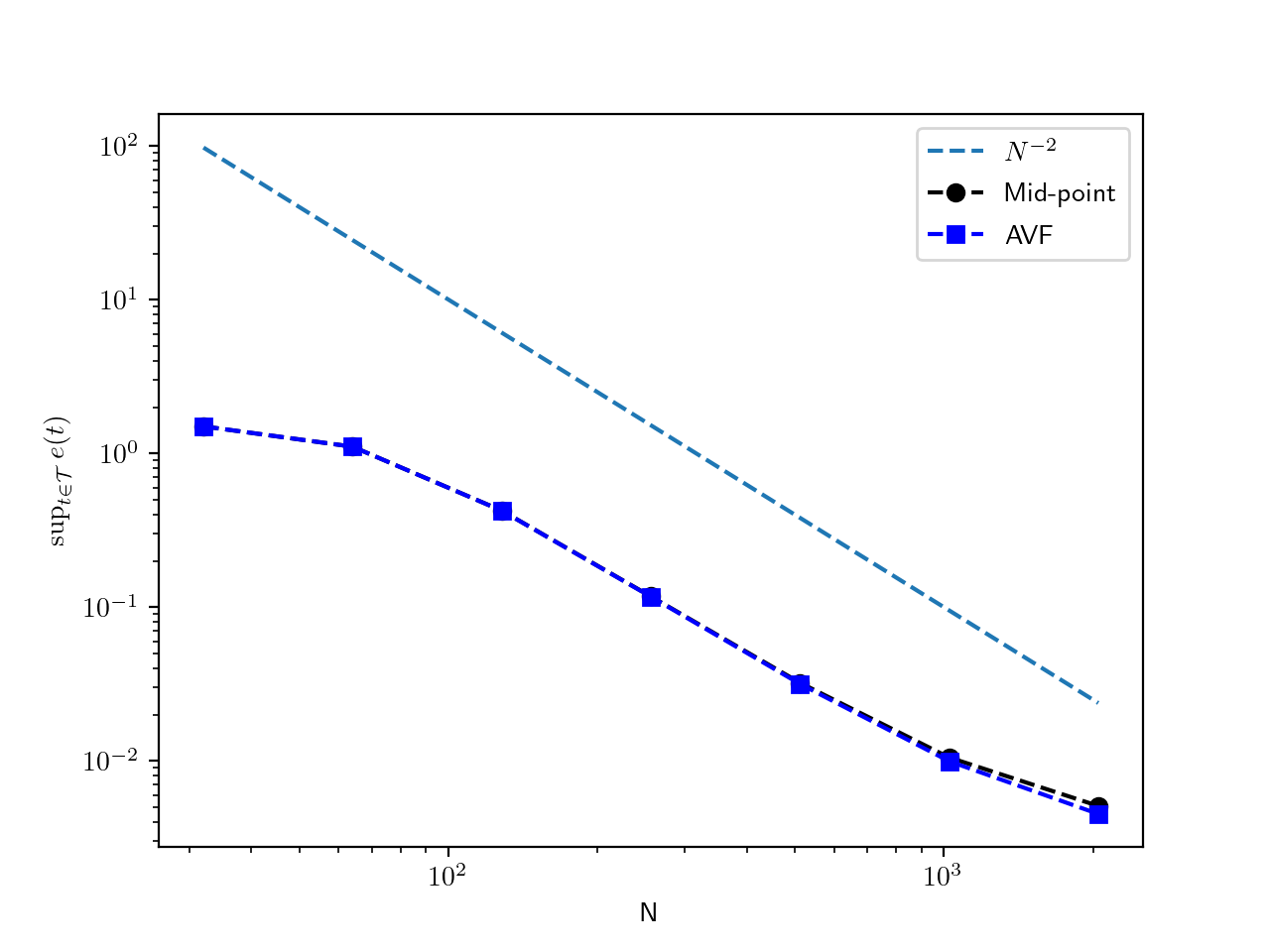}
\caption{On the left: zoom on the final time numerical solution 
compared to the analytical one. On the right, mesh convergence in log-log scale.}
\label{fig:KdV_mesh}
\end{figure}

In \Cref{fig:KdV_mesh}, on the left, we show a zoom of the solution, at final time, 
compared to the analytical one. We can see that the proposed method yields a numerical solution which does not exhibit spurious oscillations, and has a correct amplitude. The error is mainly due to a delay with respect to the analytical solution. On the same figure, on the right, we show the plot of the convergence of the solution towards the analytical one, in $L^{\infty}(\mathcal{T},L^2(\Omega))$ relative error. In this figure, we compare the errors obtained with both the time integration schemes. We clearly see that they are equivalent, the errors being slightly smaller for the AVF method. We can distinguish a pre-asymptothic phase, in which the decay of the error is slower than $N^{-2}$, the expected rate, a phase in which the schemes are second-order in space, and, lastly, for larger space resolution, we can observe that the convergence slows down. This is due to the interaction between the space discretisation and the time integration errors.  

\begin{figure}[h]
\includegraphics[scale=0.5]{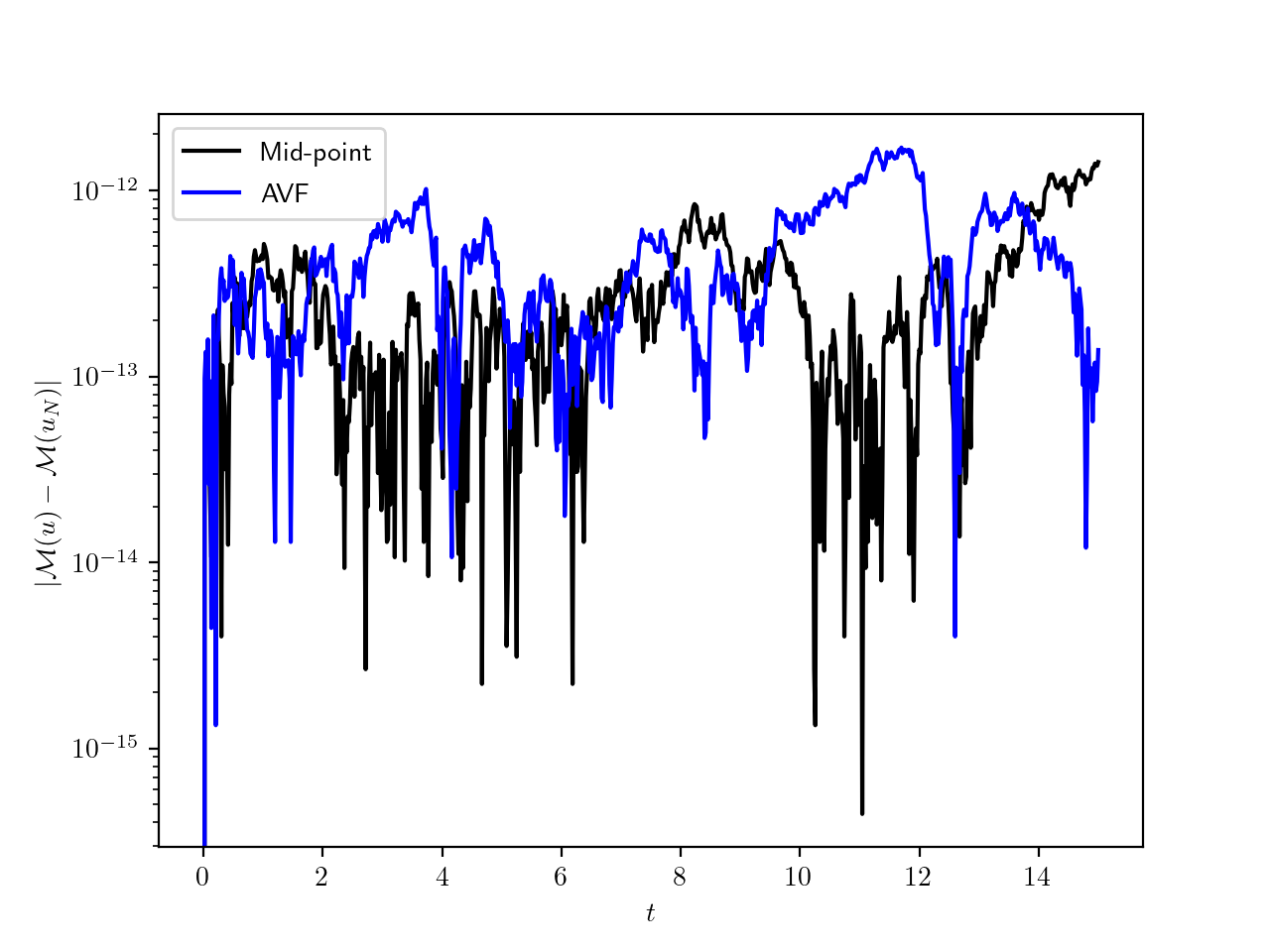}
\includegraphics[scale=0.5]{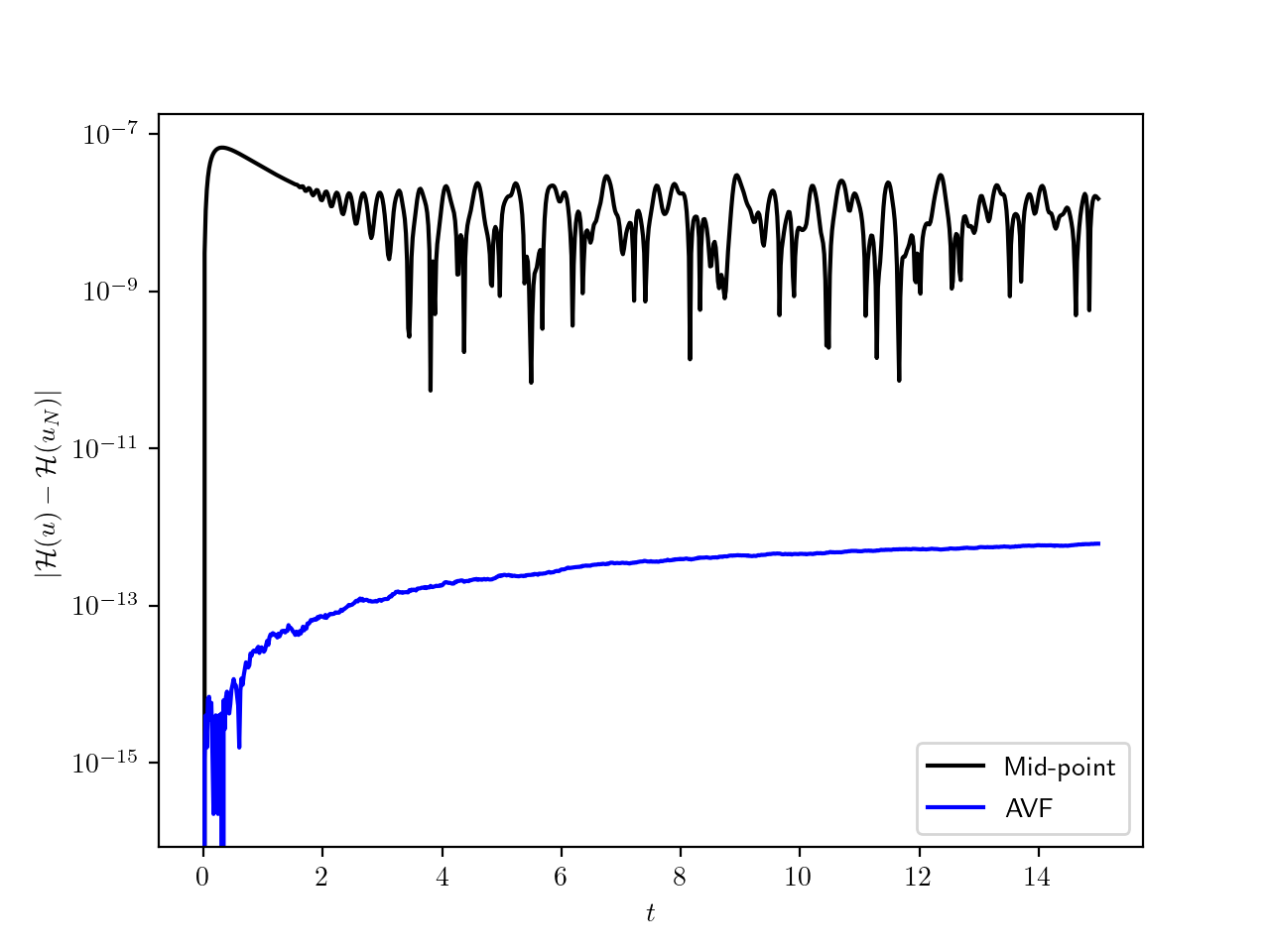}
\caption{Mass (on the left) and Hamiltonian (on the right) conservation in time, for the solution of the KdV equation with no dissipation.}
\label{fig:KdV_Ham}
\end{figure}

In \Cref{fig:KdV_Ham} we report the conservation of mass (on the left) and of the Hamiltonian (on the right).
The reference values are computed using the analytical solution \eqref{eq:soliton}.
The proposed spatial discretisation ensures that the semi-discrete system preserves the mass and the Hamiltonian. The error that it is observed in the numerical tests is purely due to the time integration scheme.
 In terms of mass conservation, both the midpoint integration method and the AVF conserve the mass up to the tolerance of the fixed-point iteration. The difference between the two time schemes can be seen on the Hamiltonian conservation, which is roughly 4 orders of magnitudes better for AVF, which conserves the Hamiltonian up to the tolerance of the nonlinear iteration. This is expected since the Hamiltonian for the KdV equation is not a quadratic function and, thus, the implicit midpoint scheme does not preserve it exactly.

\subsubsection*{Dissipative solutions}
In this section, we assess the behaviour of the method when introducing dissipation, namely by setting $\nu=1/4$. 

\begin{figure}[H]
\includegraphics[scale=0.5]{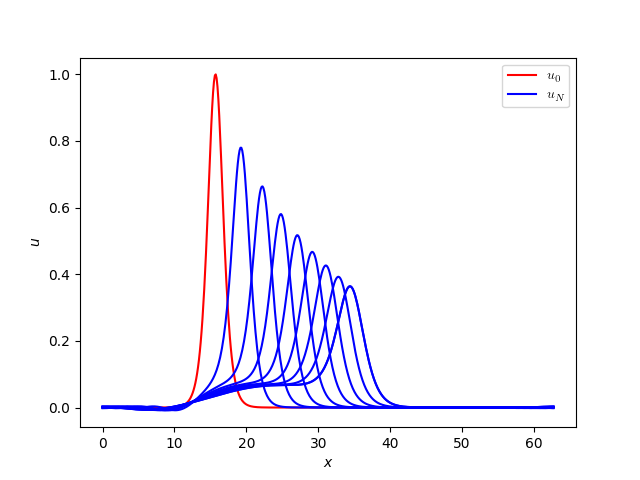}
\includegraphics[scale=0.5]{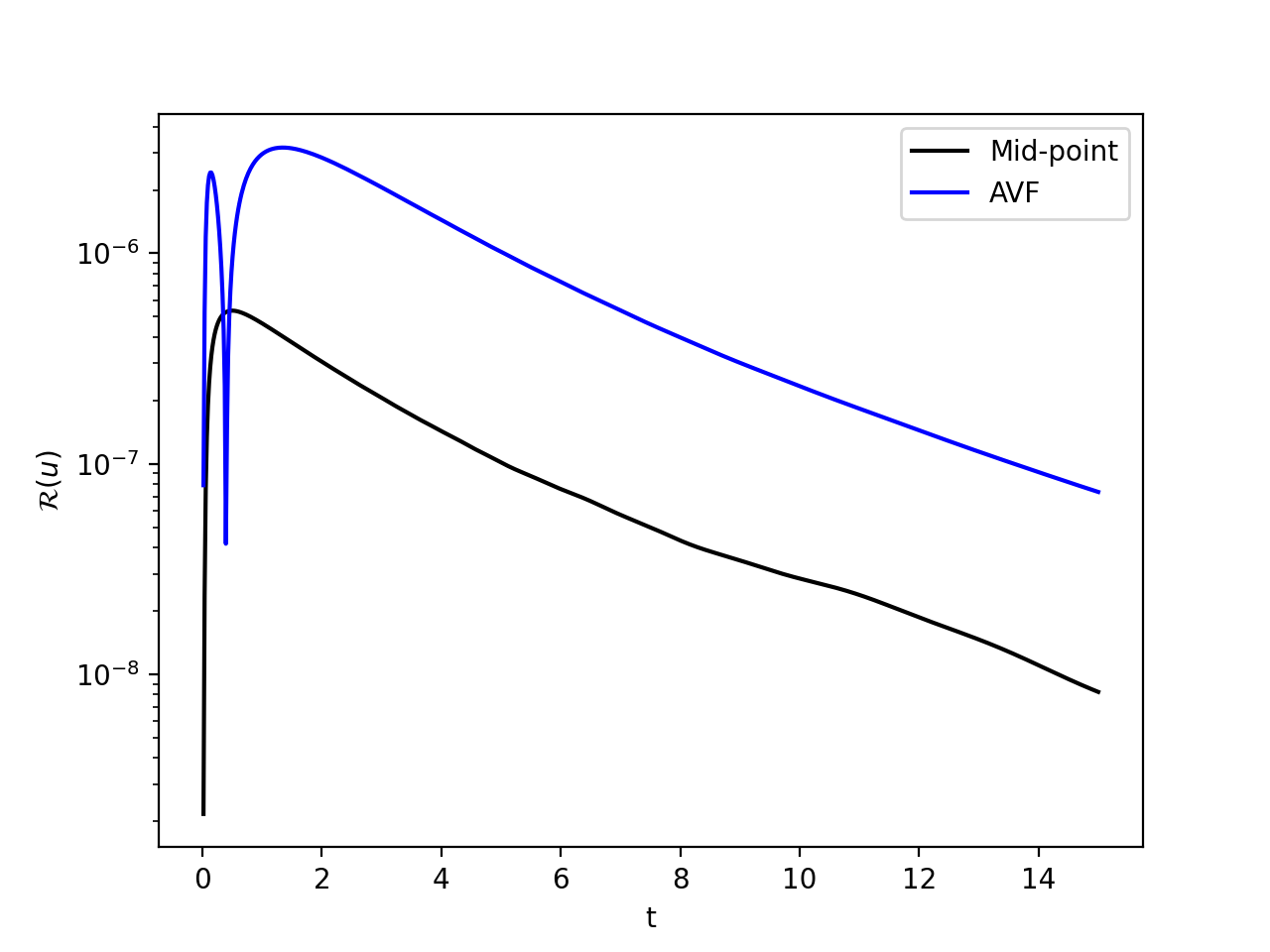}
\caption{On the left, time evolution of the solution of the dissipative KdV equation; the red curve is the initial condition evaluated at the mesh nodes.
On the right, the discrete residual in time of the entropy production rate.}
\label{fig:dissip_KdV}
\end{figure}

In \Cref{fig:dissip_KdV}, on the left, we show the time evolution of the numerical solution, obtained with 
the AVF time integration method.
We plot the solution in space every $100$ iteration, and highlight, in red, the initial condition. We can see that the dissipation mechanism induces a decrease in time of the amplitude of the soliton.
Since the proposed approximation guarantees the separation between conservative and dissipative processes, the Hamiltonian part is the classical KdV equation, and hence a decrease in amplitude is related to a decrease in the soliton speed. This phenomenon can be seen visually when considering the soliton peak at multiple times. The mass of the solution is conserved, and for the sake of brevity we do not report the plot of the mass conservation, which is analogous to the one shown for the classical KdV solution in \Cref{fig:KdV_Ham}. 

It is interesting to consider, in this case, the evolution of the entropy $\mathcal{S}(u) = -\frac{1}{2}\int_{\Omega} u(x)^2\ dx$.
We can see that, at continuous level, it holds:
$$
d_t \mathcal{S}(u) = \nu \int_{\Omega} \partial_x u(x)^2 \ dx.
$$
When considering $\nu>0$, we do not have an analytical solution for which we can exactly compute the entropy production rate. However, we can introduce the discrete residual of the entropy production
$$
\mathcal{R}(u_N^{k+1}) := \frac{\Scal(u_N^{k+1}) - \Scal(u_N^{k})}{\Delta t} - \nu \int_{\Omega} \left(\partial_x \Big(\frac{u_N^{k+1}+u_N^{k}}{2} \Big)\right)^2 \ dx.
$$
The evolution of this quantity in time is shown, for both time integration schemes, in \Cref{fig:dissip_KdV}, on the right. We can see that both methods have a good agreement, at discrete level, in terms of entropy production. 


\subsection{The 2D Navier-Stokes equations in vorticity formulation}
We recall the formulation of the 2D Navier-Stokes equations as an infinite dimensional Hamiltonian flow, the incompressible Euler equations, to which we add a dissipation mechanism. 

Let $\Omega\subset\mathbb{R}^3$ be a two-dimensional spatial domain.
Let $\omega:\Tcal\times \Omega\rightarrow\mathbb{R}$ be the vorticity defined as
$\omega=\curl u= \DIV(\mathbb{J} u)$ where
$u$ is the velocity field of the flow and
$$\mathbb{J}:=
\begin{pmatrix}
    0 & \mathrm{Id}\\
    -\mathrm{Id} & 0
\end{pmatrix}.$$
Introducing the stream function $\Psi:\Tcal\times \Omega\rightarrow\mathbb{R}$ defined as $\curl \Psi=-\mathbb{J}\nabla\Psi= u$, then the 2D incompressible Navier-Stokes equations read
\begin{equation*}
    \begin{cases}
        \partial_t\omega + \{\Psi,\omega\} = \nu\Delta\omega\\
        \omega=\Delta\Psi
    \end{cases}
\end{equation*}
where $\nu$ is a positive constant, the stream function $\Psi$ is related to the vorticity via the Laplace-Beltrami operator $\Delta$, and 
$\bra{\cdot}{\cdot}$ is the Poisson bracket defined as
$\bra{\Psi}{\omega}=\nabla\Psi\cdot\mathbb{J}\nabla\omega$. 

The system can be written as a Hamiltonian system with dissipation; thereby
\begin{equation}\label{eq:NSHamDiss}
    \partial_t\omega=\Jc(\omega)\dH(\omega)+\nu \Gc\delta\Ecal(\omega),
\end{equation}
where the function $\Hcal$ is the kinetic energy
$$
\Hcal(\omega)=-\dfrac12\int_{\Omega}\Psi\omega\, dx,
\qquad\delta\Hcal=-\Psi
$$
and $\Ecal$ is the enstrophy
$$\Ecal(\omega)=\dfrac12\int_{\Omega}|\omega|^2\, dx,\qquad\delta\Ecal=\omega.$$
The metric tensor $\Gc=\Delta$ is the Laplacian, while the Poisson tensor is defined as $\Jc(\omega) z=\{z,\omega\}=\mathrm{ad}^*_{z}\omega=-\nabla\omega\cdot\curl z=-\nabla\omega\cdot\mathbb{J}\nabla z$.
Observe that the enstrophy $\Ecal$ satisfies
$$\Jc(\omega)\delta\Ecal(\omega)=\Jc(\omega)\omega=\{\omega,\omega\}=0,$$
and it is, thus, a Casimir of the $\Jc$ operator. Moreover, it dissipates since
$$d_t\Ecal(\omega)=\langle\delta\Ecal(\omega),\Jc(\omega)\dH(\omega)+\nu\Delta\delta\Ecal(\omega)\rangle_0=-\nu\int_{\Omega}|\nabla\omega|^2\, dx=-2\nu \Pcal(\omega(t)),$$
where $\Pcal$ is called palinstrophy.
The kinetic energy $\Hcal$ satisfies
$$\Gc\dH(\omega)+\omega=-\Delta\Psi+\omega=0$$
and it dissipates since
$$d_t\Hcal(\omega)
=\langle\dH(\omega),\nu\Delta\delta\Ecal(\omega)\rangle_0
=-\nu\int_{\Omega}|\omega|^2\, dx=-\nu 2 \Ecal(\omega(t)).$$

Let $X_s=X_1=H^1(\Omega)$ and $X_{s-d}=X_{s-e}=X_1$.
Following \Cref{sec:mixed}, the mixed formulation of problem \eqref{eq:NSHamDiss} reads: for any $v\in X_1$,
\begin{equation*}
\begin{cases}
\langle \partial_t \omega, v \rangle_0 = -\langle z, \Jc(\omega) v \rangle_0 - \nu \langle \nabla y, \nabla v \rangle_0,\\
\langle z, v\rangle_0 = -\langle \Psi, v \rangle_0,\\
\langle y, v \rangle_0 = \langle \omega, v \rangle_0,
\end{cases}
\end{equation*}
where the stream function is solution of the Laplace equation $\inpo{\omega,v}=-\inpo{\nabla\Psi,\nabla v}$ for any $v\in X_1$.
From this formulation, we can check that the choice $V_N = \mathbb{P}^1(\Omega_N)$ is a simple pertinent choice to discretise the system. We will use classical Lagrange piecewise linear finite element to discretise both the vorticity and the stream function. 

Using the notation from \Cref{sec:confdisc}, the semi-discrete 
system reads
\begin{equation}\label{eq:NSsemi-disc}
    \begin{cases}
        M d_t a=\Jd(a) b -\nu K a\\
        K b = M a
    \end{cases}
\end{equation}
where $M, K\in\mathbb{R}^{N\times N}$ are the mass and stiffness matrix, respectively, while
$$\Jd(a)_{i,j}=\sum_{k=1}^N a_k\inpo{\{v_k,v_j\},v_i}
=\sum_{k=1}^N a_k \inpo{\nabla v_k\cdot \mathbb{J}\nabla v_j,v_i}\qquad 1\leq i,j\leq N.$$

We verify that \eqref{eq:NSsemi-disc} has the Hamiltonian dissipative structure and preserves the conservation and dissipation laws at the semi-discrete level.
First, observe that the discrete Hamiltonian and discrete entropy are
\begin{equation}\label{eq:discrHE}
    \begin{aligned}
        & H(a)=-\dfrac12 \langle\omega_N,\Psi_N\rangle_0=\dfrac12 a^{\top}M b=\dfrac12 a^{\top}MK^{-1}Ma, & \qquad\nabla H(a)&= MK^{-1}Ma, \\
        & E(a) = \dfrac12 a^{\top} Ma, & \qquad\nabla E(a) &= Ma.
    \end{aligned}
\end{equation}
With these definitions, equation \eqref{eq:NSsemi-disc} can be written in the form
\begin{equation*}
    d_t a = \Jm(a)\nabla H(a)+\nu \Gm\nabla E(a),
\end{equation*}
where $H$ and $E$ are as in \eqref{eq:discrHE}, $\Gm=-M^{-1}KM^{-1}$ is symmetric negative definite, 
$\Jm(a)$
is defined as the matrix-valued function $\Jm:a\in\mathbb{R}^N\mapsto \Jm(a)\in\mathbb{R}^{N\times N}$ such that
$\Jm(a)=M^{-1}\Jd(a)M^{-1}$ or componentwise
$$\Jm(a)_{i,j}:=\sum_{k,\ell,s=1}^N  M^{-1}_{i,s} M^{-1}_{\ell, j} a_k \inpo{\{v_k,v_{\ell}\},v_s}\qquad 1\leq i,j\leq N.$$
The matrix $\Jm(a)$ is skew-symmetric for any $a\in\mathbb{R}^N$
and it satisfies the Jacobi identity, see \Cref{lem:poisson_bracket}.

We check the null condition, namely $\Jm(a)\nabla E(a)=0$ for any $a\in\mathbb{R}^N$,
\begin{equation*}
    \begin{aligned}
        (\Jm(a)Ma)_i & 
        =  \sum_{k,s=1}^N M^{-1}_{i,s}a_k\Jd(a)_{s,k}
         = \sum_{k,\ell,s=1}^N M^{-1}_{i,s} a_ka_{\ell} 
        \inpo{\{v_{\ell},v_k\}, v_s}\\
         & = \sum_{s=1}^NM^{-1}_{i,s} \langle\{\omega_N,\omega_N\}, v_s\rangle_0 = 0.
    \end{aligned}
\end{equation*}
The other condition is
$$\Gm\nabla H(a)+a=-M^{-1}KM^{-1}(MK^{-1}Ma)+a=0.$$

For the numerical time integration of the semi-discrete problem \eqref{eq:NSsemi-disc}, we consider the implicit midpoint rule and the AVF method described in \Cref{sec:time}; thereby
\begin{equation*}
    \begin{aligned}
        \left(M+\dfrac{\dt}{2}\nu K\right) a^{k+1}  
        & -d_1^{-1} \dt\Jd(a^{k+1})b^{k+1}
        -d_2^{-1}\dt\Jd(a^{k+1})b^k
        -d_2^{-1}\dt \Jd(a^{k})b^{k+1}& \\
        & =\left(M-\dfrac{\dt}{2}\nu K\right) a^{k}
+d_1^{-1}\dt \Jd(a^{k})b^{k}.
    \end{aligned}
\end{equation*}
The implicit midpoint rule corresponds to $d_1=d_2=4$; whereas for the AVF scheme to $d_1=3$, $d_2=6$.

Since the Hamiltonian \eqref{eq:discrHE} of the Navier-Stokes problem is a quadratic function, in the absence of dissipation, it is conserved by the implicit midpoint rule as shown in \Cref{lem:IMconsHam}.

\subsubsection*{The Navier-Stokes equations on the torus}
When using periodic boundary conditions, it is possible to have analytical solutions of the Navier-Stokes equations, similar to the classical Taylor vortex solutions. We refer to \cite{Walsh92} for further details.

Let $\Omega=[0,2\pi]^2$ and $V=H^1_{\per}(\Omega)$ and let the temporal interval be $\Tcal=(0,1]$. We take $\nu\in\mathbb{R}^+$ and $\lambda=25$. We consider the stream function given by
$$
\Psi(t,x,y) = e^{-\nu\lambda t} \left[ \frac{1}{4}\cos(3x)\sin(4y) - \frac{1}{5}\cos(5y) - \frac{1}{5}\sin(5x) \right],\qquad (t,x,y)\in\Tcal\times \Omega.
$$
We can verify that the vorticity $\omega$ satisfies
$$
\omega = -\lambda \Psi.
$$
The evolution is such that the velocity field is always orthogonal to the gradient of the vorticity, and hence we can immediately verify that the above-stated solution verifies the Navier-Stokes equation in vorticity formulation. 

In the present case, we consider several values of the diffusivity, $\nu\in\left\lbrace 10^{-2}, 10^{-4}, 0 \right\rbrace$. The last value corresponds to the steady solution of the incompressible Euler flow, for which $u\cdot \nabla \omega = 0$, and $\partial_t\omega = 0$. This limit case is challenging as there are no dissipation mechanisms acting.
We consider $N_t=200$ equispaced time steps.

\begin{figure}[H]
\includegraphics[scale=0.55]{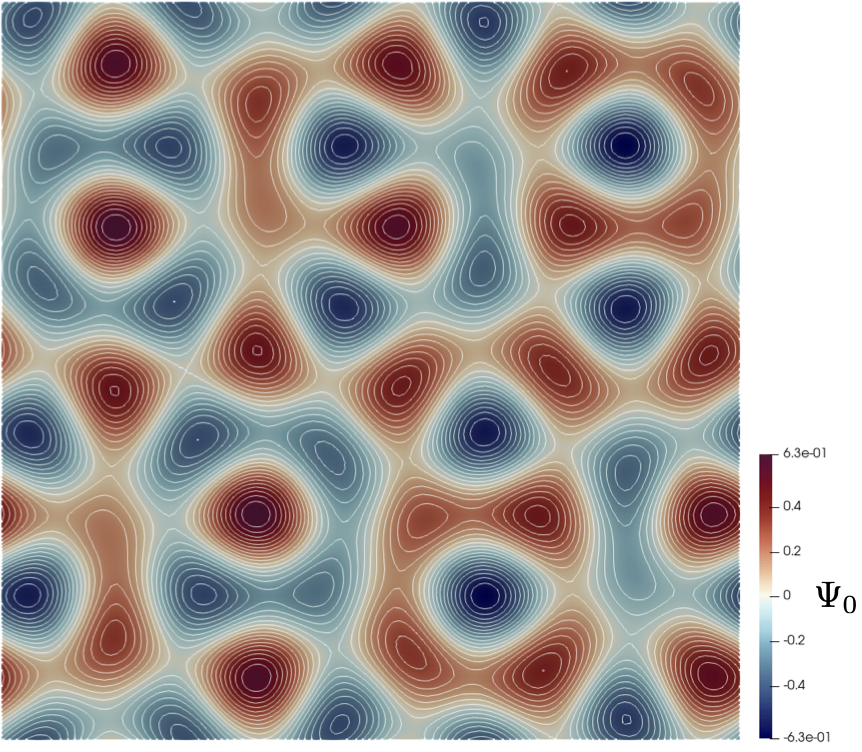}
\includegraphics[scale=0.55]{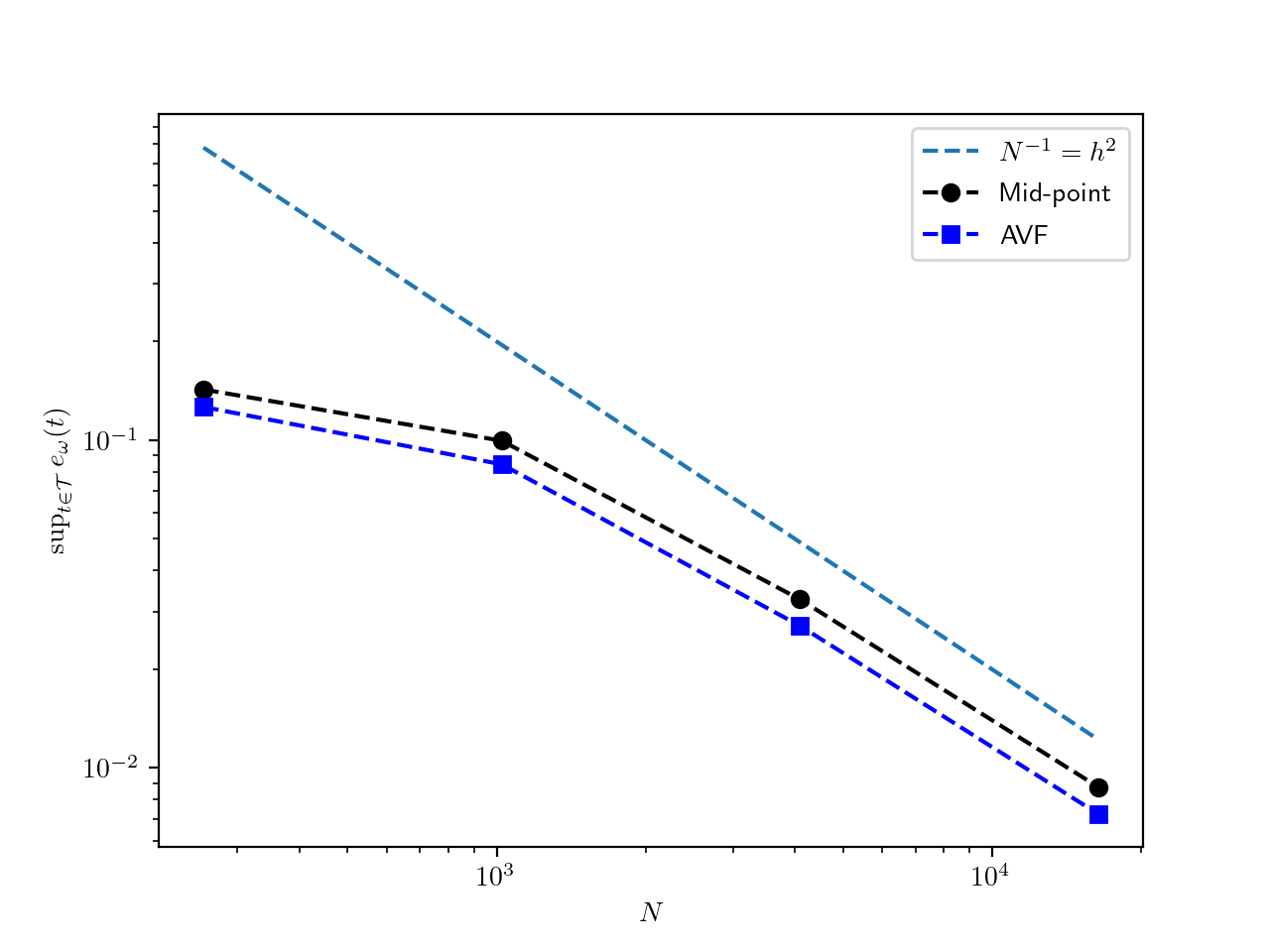}
\caption{Initial condition for the stream function $\Psi$ (on the left), and convergence in mesh (on the right), for the Navier-Stokes equation solution on the torus, when $\nu=10^{-2}$.}
\label{fig:walsh_1}
\end{figure}
In \Cref{fig:walsh_1}, on the left, we show the initial condition for the stream function $\Psi_0$. As we can see, it is periodic, but it does not have the symmetries of the classical Taylor vortex solution. On the right, the convergence in mesh is shown. In particular, we plot the error on the vorticity in $L^{\infty}(\mathcal{T}, L^2(\Omega))$ when $\nu=10^{-2}$. After a pre-asymptotic phase, we can see that we recover the expected rate $h^2=N^{-1}$. The results are analogous when taking $\nu=10^{-4}$. 

\begin{figure}[H]
\includegraphics[scale=0.5]{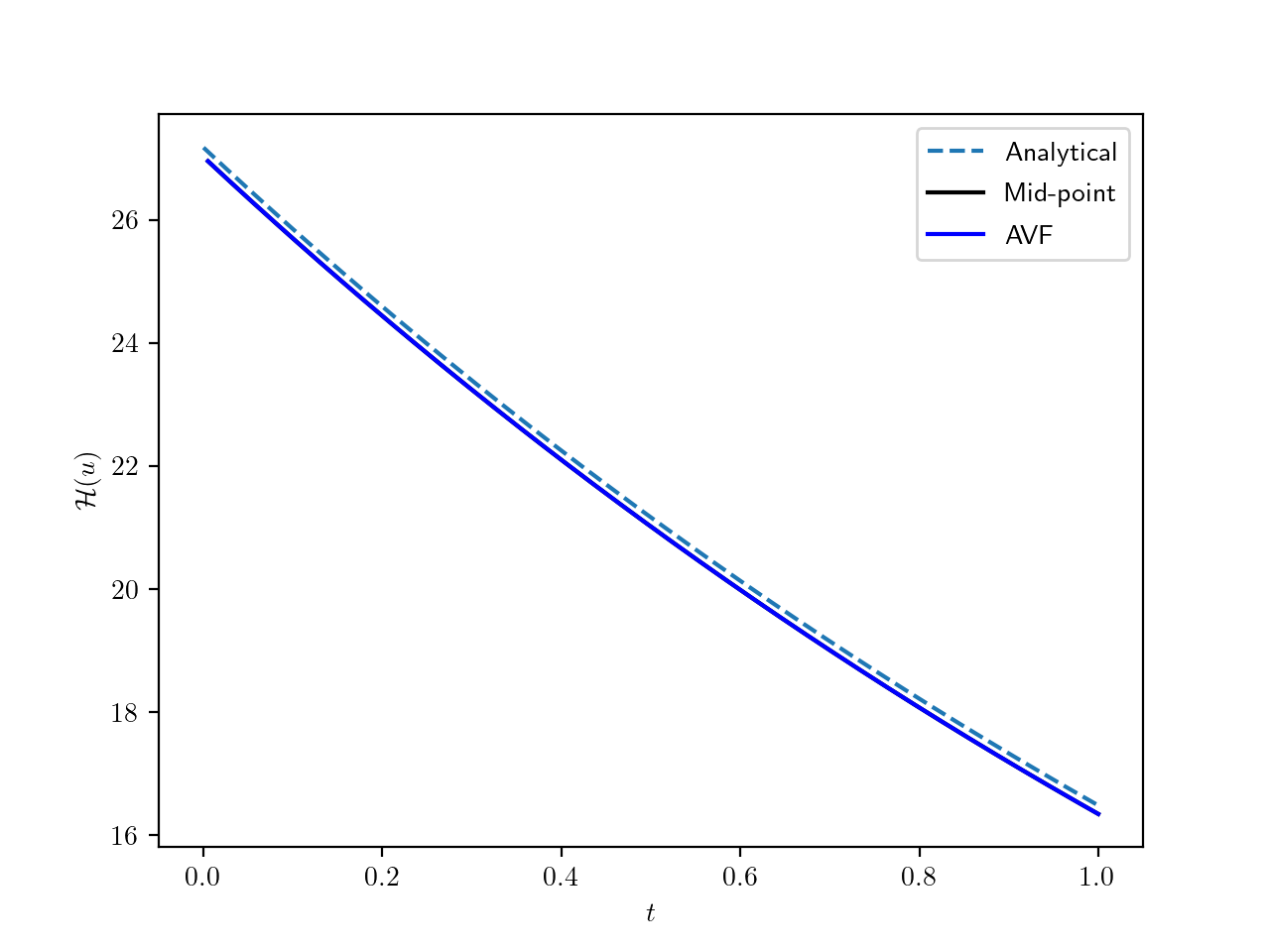}
\includegraphics[scale=0.5]{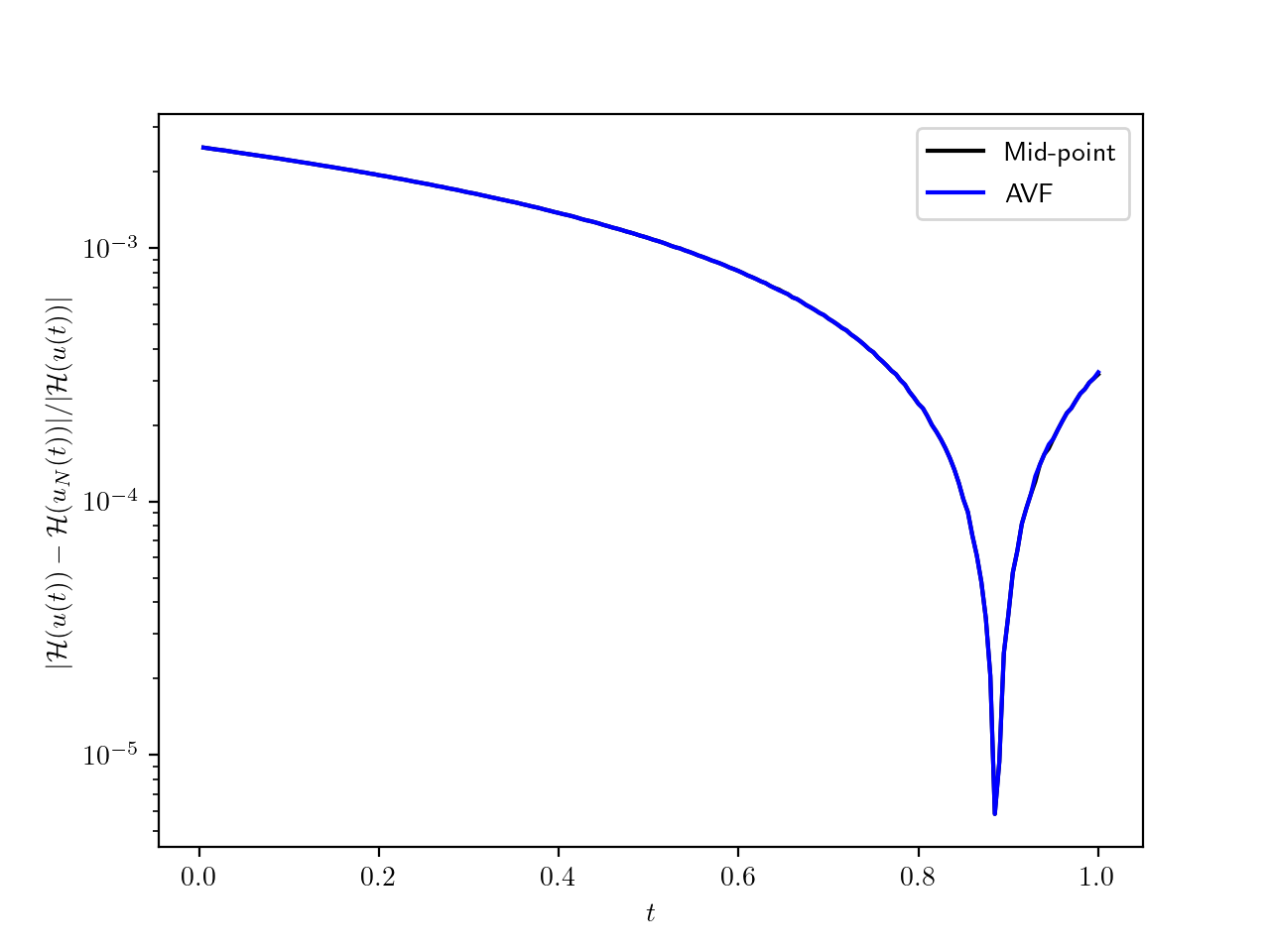}
\caption{Hamiltonian decay (left), and error in semi-logartithmic scale (right).}
\label{fig:Ham_NS_1em2}
\end{figure}
In \Cref{fig:Ham_NS_1em2}, on the left, we show the decay of the Hamiltonian of the flow when $\nu=10^{-2}$, and $N=128^2$. We can see that both the AVF scheme and the implicit midpoint rule have the same behaviour, and follow the analytical trend. The error we have on the analytical value of the Hamiltonian depends solely on the fact that the initial condition does not belong to the space $V_N$, and hence we have a discretisation error. On the right of \Cref{fig:Ham_NS_1em2}, we show the evolution of the error on the Hamiltonian decay, that is the difference between the values of the Hamiltonian $\mathcal{H}(\omega_N)$ obtained by the numerical scheme and the theoretical decay that the Hamiltonian should have, namely $\mathcal{H}(\omega_N(0))e^{-2\lambda t}$. The plot shows that both schemes have a good agreement (although not to machine precision) with respect to the theoretical decay. 

\begin{figure}[H]
\includegraphics[scale=0.5]{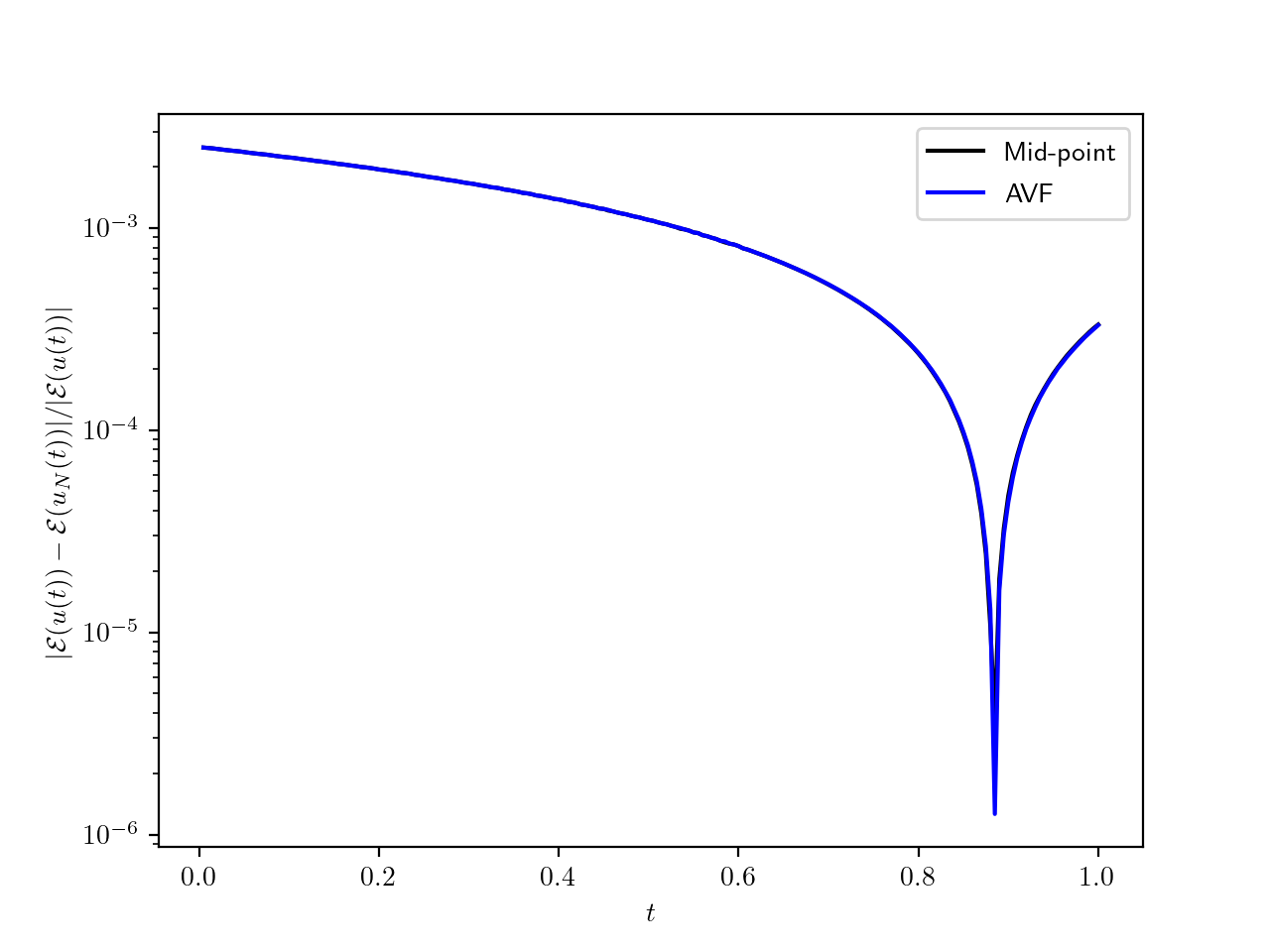}
\includegraphics[scale=0.5]{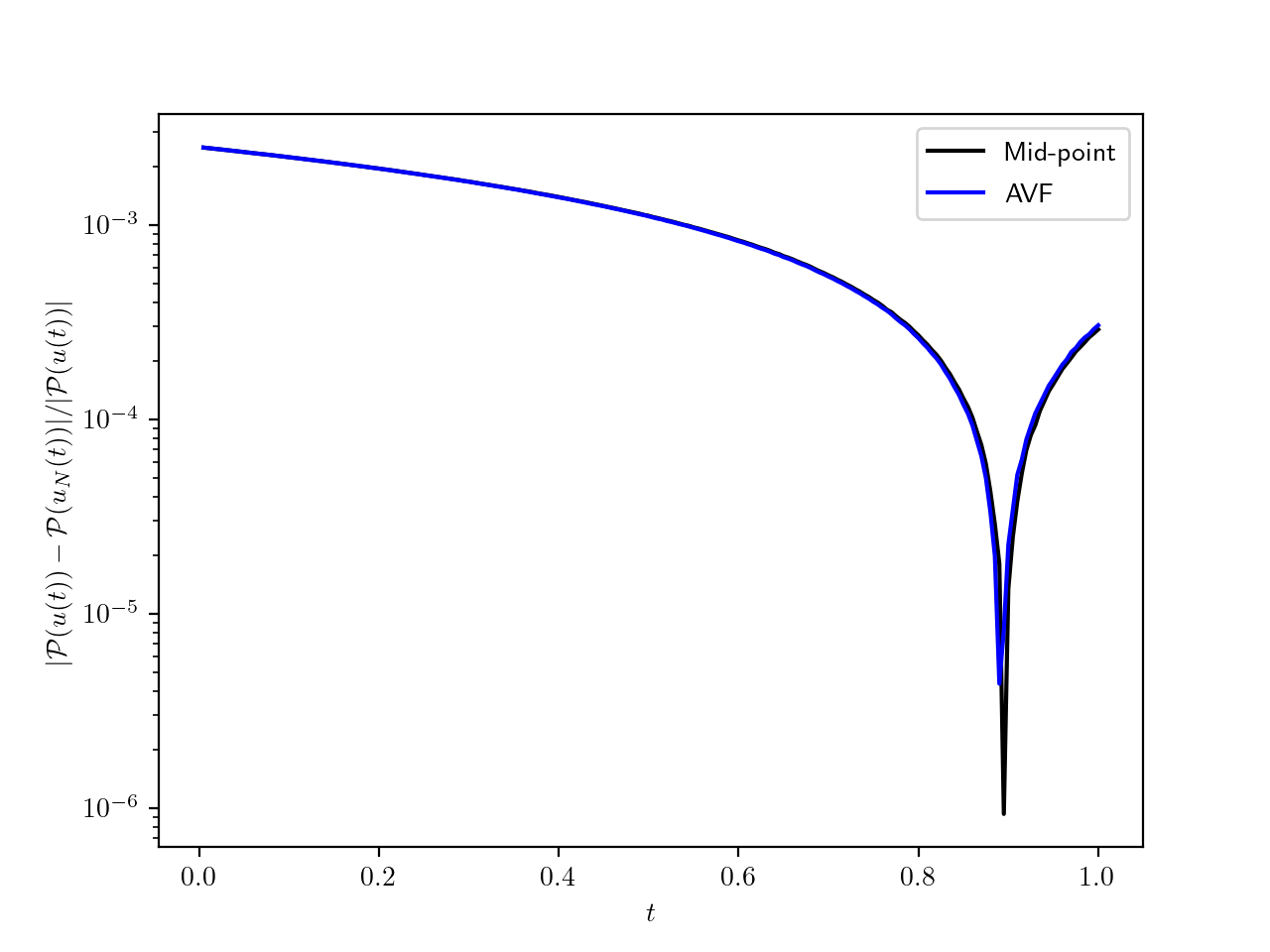}
\caption{Evolution of the error in the enstrophy (left), and in the palinstrophy (right).}
\label{fig:EnPal_NS_1em2}
\end{figure}
In \Cref{fig:EnPal_NS_1em2} we show the same plot for the error in the decay of enstrophy and palinstrophy. For both these quantities, whose analytical decay in time is the same as for the Hamiltonian, we can see the same behaviour, for both time integration schemes. 

\begin{figure}[H]
\includegraphics[scale=0.5]{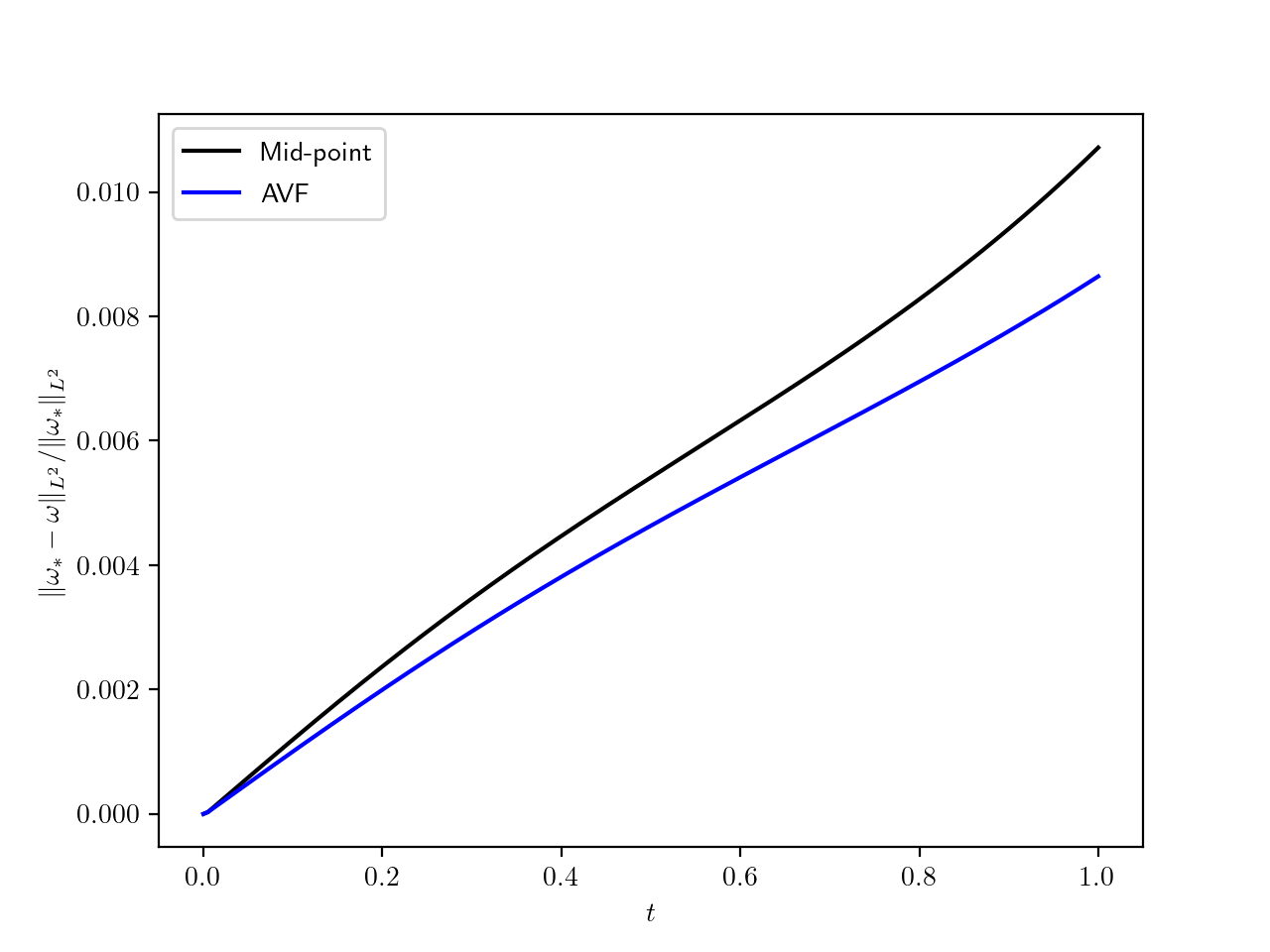}
\includegraphics[scale=0.5]{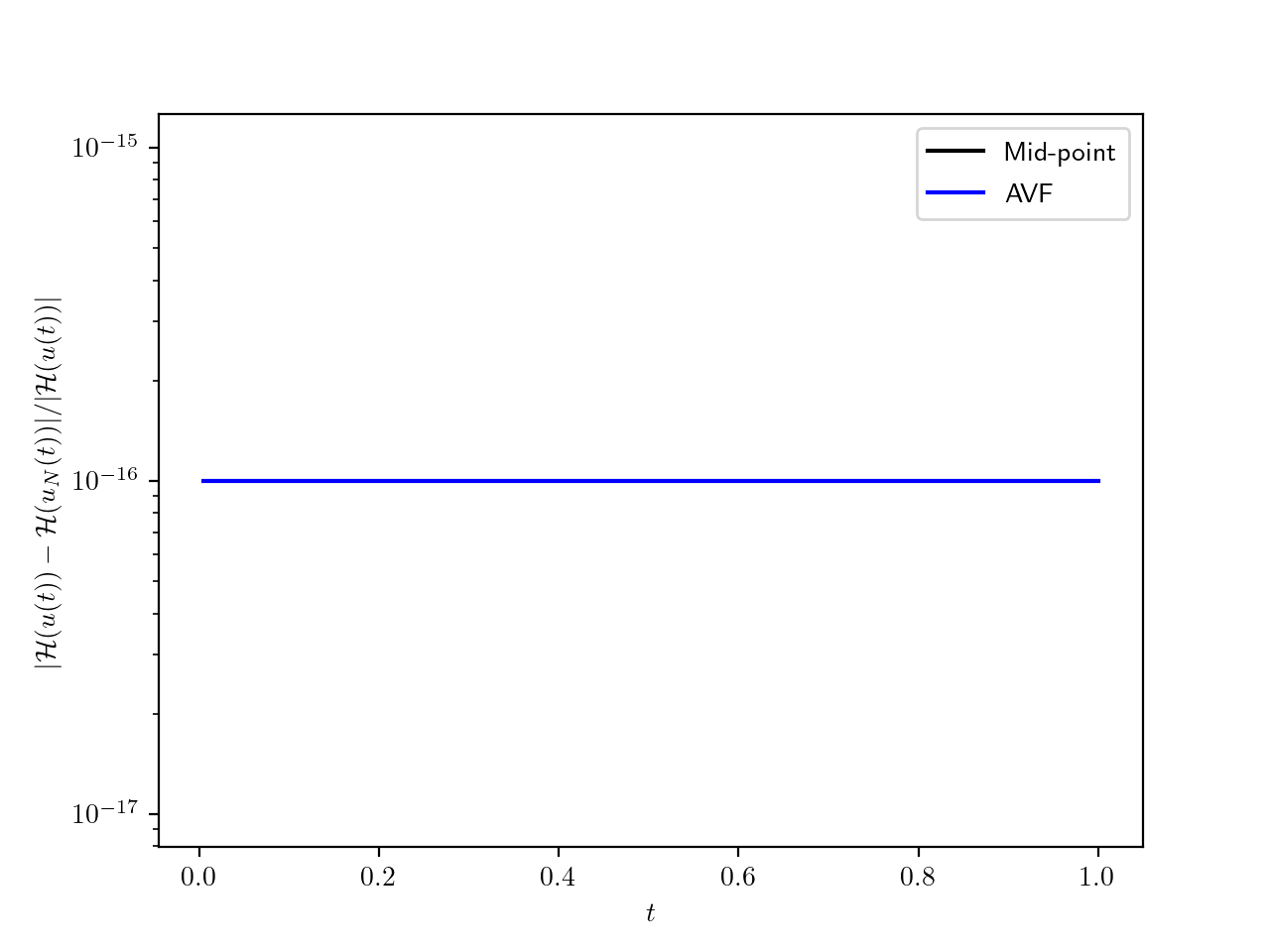}
\caption{On the left, evolution of the relative error in the approximation of the vorticity. On the right, error in the Hamiltonian evolution.}
\label{fig:err_NS_noDissip}
\end{figure}

Let us investigate the behaviour of the method when taking the limit case $\nu=0$. 
In \Cref{fig:err_NS_noDissip}, on the left, we show the evolution of the relative error in the vorticity approximation, as function of time, when taking $N=128^2$. We can see that the maximum in time of the relative error is similar to what we got in the dissipative case. On the left, we plot the error with respect to the analytical evolution of the Hamiltonian, which is constant in time since there is no dissipation. We can see that this error is of the order of the machine precision, for both the time integration schemes. The method is able to correctly represent the fact that the Hamiltonian is a constant of motion. The same happens for the enstrophy. The error on the solution is mainly related to the fact that the initial condition, which is an equilibrium point for the infinite dimensional Hamiltonian system does not belong to $V_N$. The projection on $V_N$ of the initial condition is not an equilibrium point for the discrete Hamiltonian system, and this is why the solution, and hence the error, evolves in time.  

In conclusion, the behaviour of the semi-discretisation in space is rather uniform with respect to the level of dissipation which we consider, including the case in which we do not have any dissipation.  

\subsubsection*{The Navier-Stokes equations on the sphere}
In this section, we consider the Navier-Stokes equations on the sphere. The equations are the same, with the only difference that the gradient is replaced by the covariant gradient and the Laplace operator by the Laplace-Beltrami operator on the sphere. 

Let the coordinates on the sphere be the angles $(\vartheta,\varphi)\in [0,\pi]\times[0,2\pi]$. As an example of analytical solution, we consider the stream function at initial time:
$$
\Psi_0(\vartheta,\varphi) = \frac{1}{2} \sin(\vartheta)\cos(\varphi).
$$
This function is related to spherical harmonics, and we can prove that we obtain a solution of the form:
$$
\Psi(t,\vartheta,\varphi) = \Psi_0 e^{-2 \nu t},
$$
and $\omega(t,\vartheta,\varphi) = -2 \Psi(t,\vartheta,\varphi)$.

\begin{figure}[H]
\includegraphics[scale=0.45]{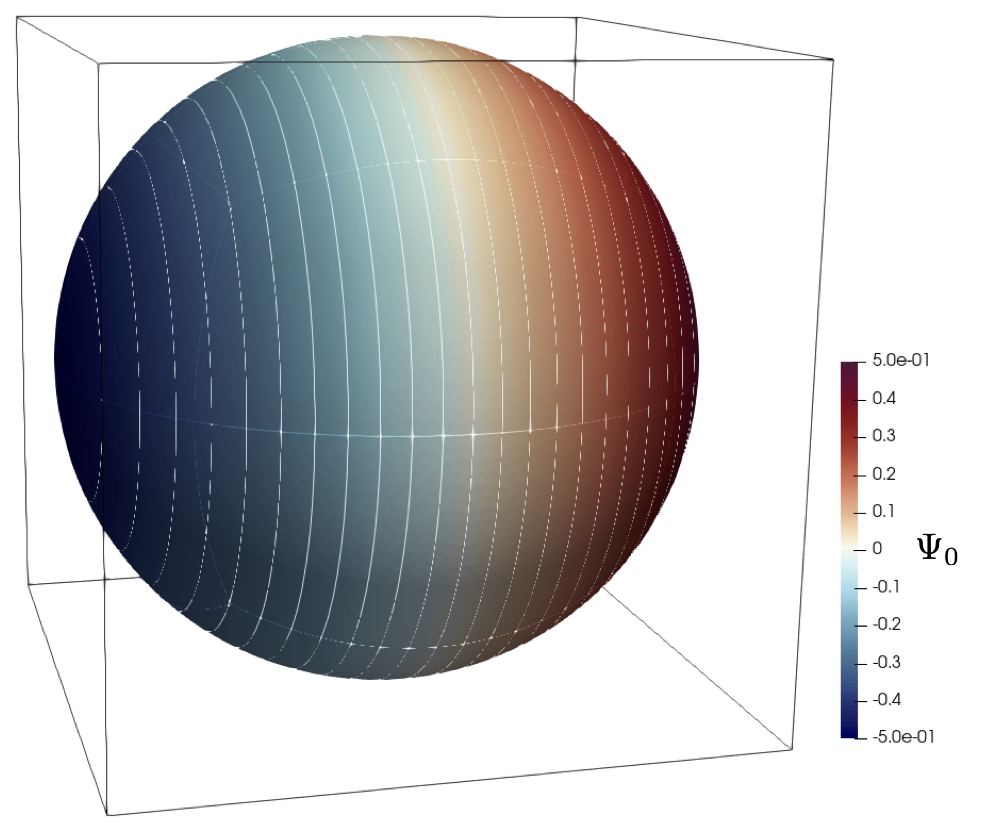}
\includegraphics[scale=0.6]{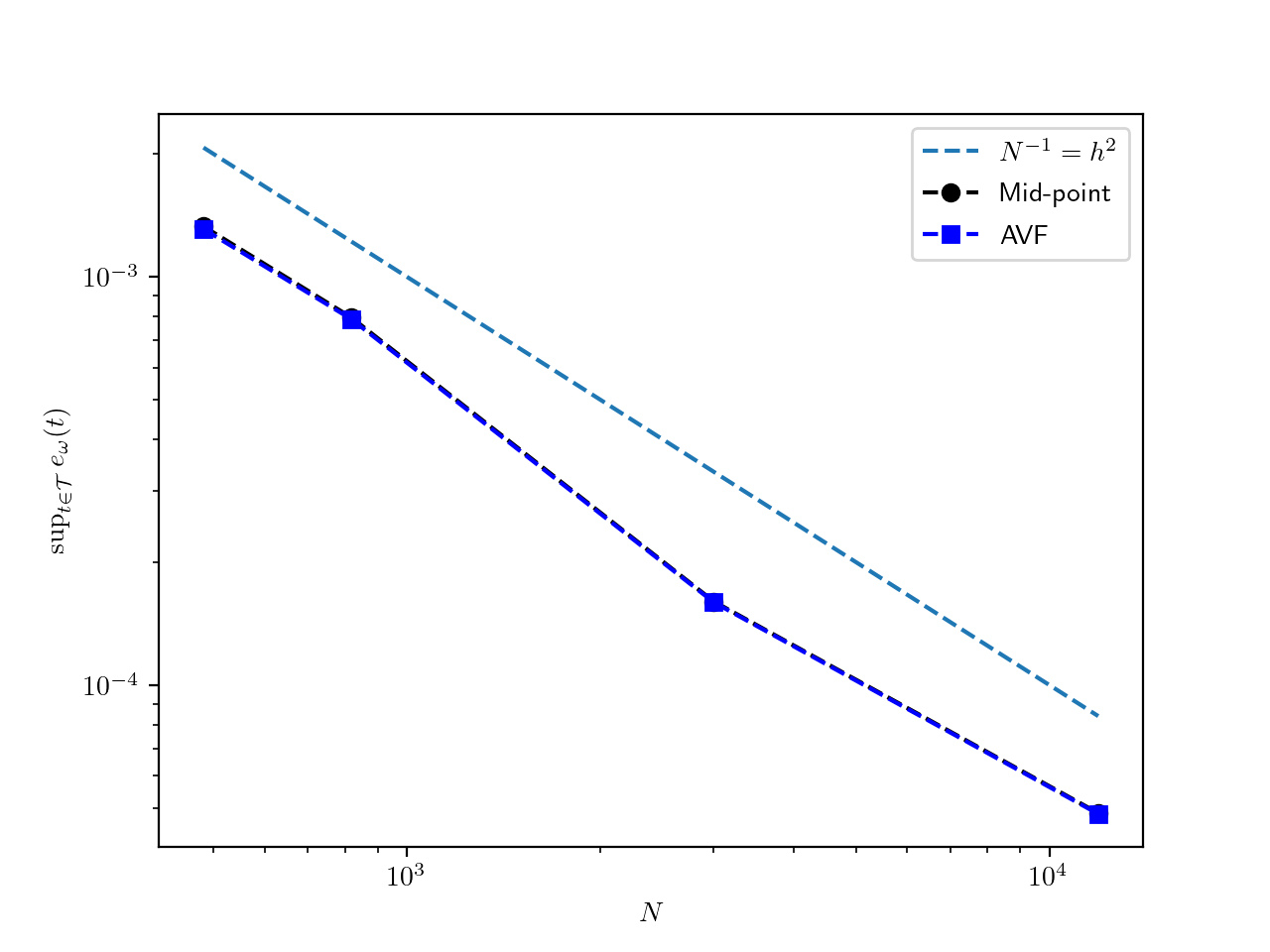}
\caption{On the left, initial condition for the stream function. On the right, convergence in mesh, for the two time integration schemes, when varying the space resolution.}
\label{fig:sphere_1}
\end{figure}
The numerical simulations are performed on a time interval $\mathcal{T}=[0,1]$, by taking $\Delta t = 10^{-2}$. The dissipation coefficient is taken, as before $\nu \in\left\lbrace 10^{-2}, 10^{-4}, 0 \right\rbrace$. Both the implicit midpoint rule and the AVF time integration methods are tested.

In \Cref{fig:sphere_1}, on the left, we show the function $\Psi_0$ on the sphere. On the right, we report the plot of the relative $L^{\infty}(\mathcal{T}, L^2(\Omega))$ error as function of the number of degrees of freedom used to discretise the mesh. As we can see, we obtain the same behaviour as for the planar case considered in the previous section. There are no substantial differences concerning the time integration scheme. The errors plotted in the figure are the ones obtained for $\nu=10^{-2}$. The behaviour of the error obtained for $\nu=10^{-4}$ and in the inviscid case are analogous and are not reported for the sake of brevity. 

As a second test case, we are going to consider a scenario in which we do not have an analytical solution. 
\begin{figure}[H]
\includegraphics[scale=0.6]{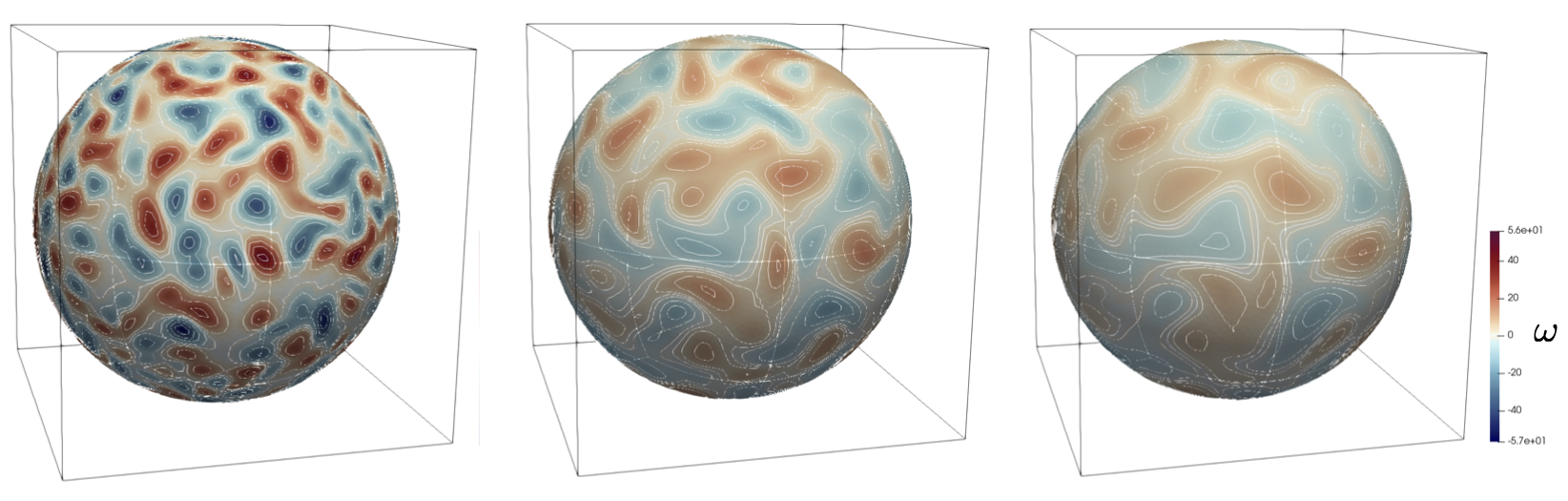}
\caption{Evolution in time of the vorticity}
\label{fig:vortices}
\end{figure}
In particular, we initialise the vorticity as a sum of $512$ point vortices, half of which with positive sign, the other half with negative sign, and intensity $\overline{\omega}_i = 400$. 

In \Cref{fig:vortices} we show the solution obtained on the finest mesh $N\approx 1.2\cdot 10^4$, when $\nu = 10^{-2}$, at $t=\left\lbrace 0.1,0.3,0.5 \right\rbrace$ using the midpoint rule as time integration scheme. In this figure we can see the dissipation of the enstrophy as well as the vortices coalescence phenomenon which is typical of two dimensional flows. 

\begin{figure}[H]
\includegraphics[scale=0.35]{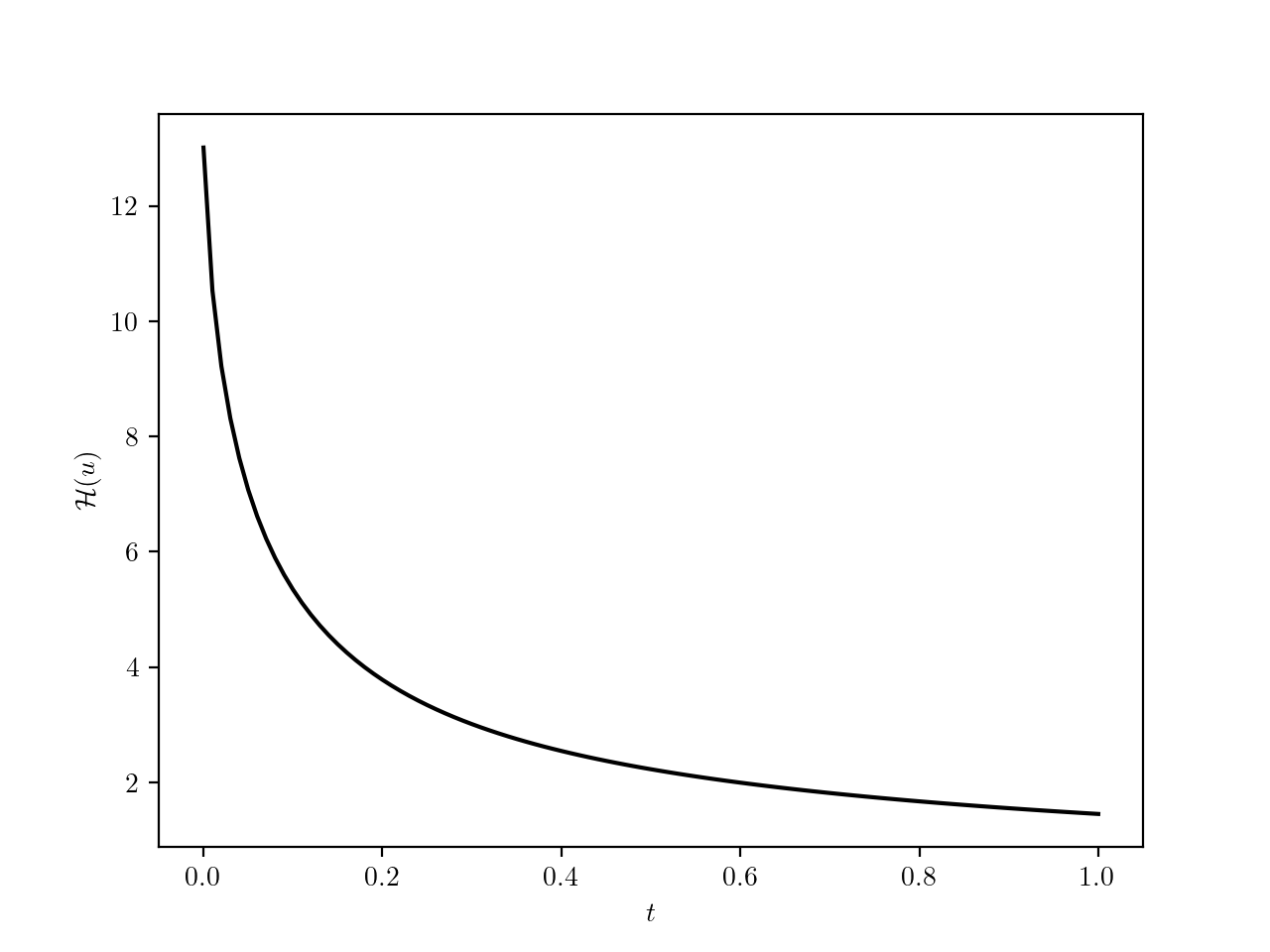}\hspace{-2ex}
\includegraphics[scale=0.35]{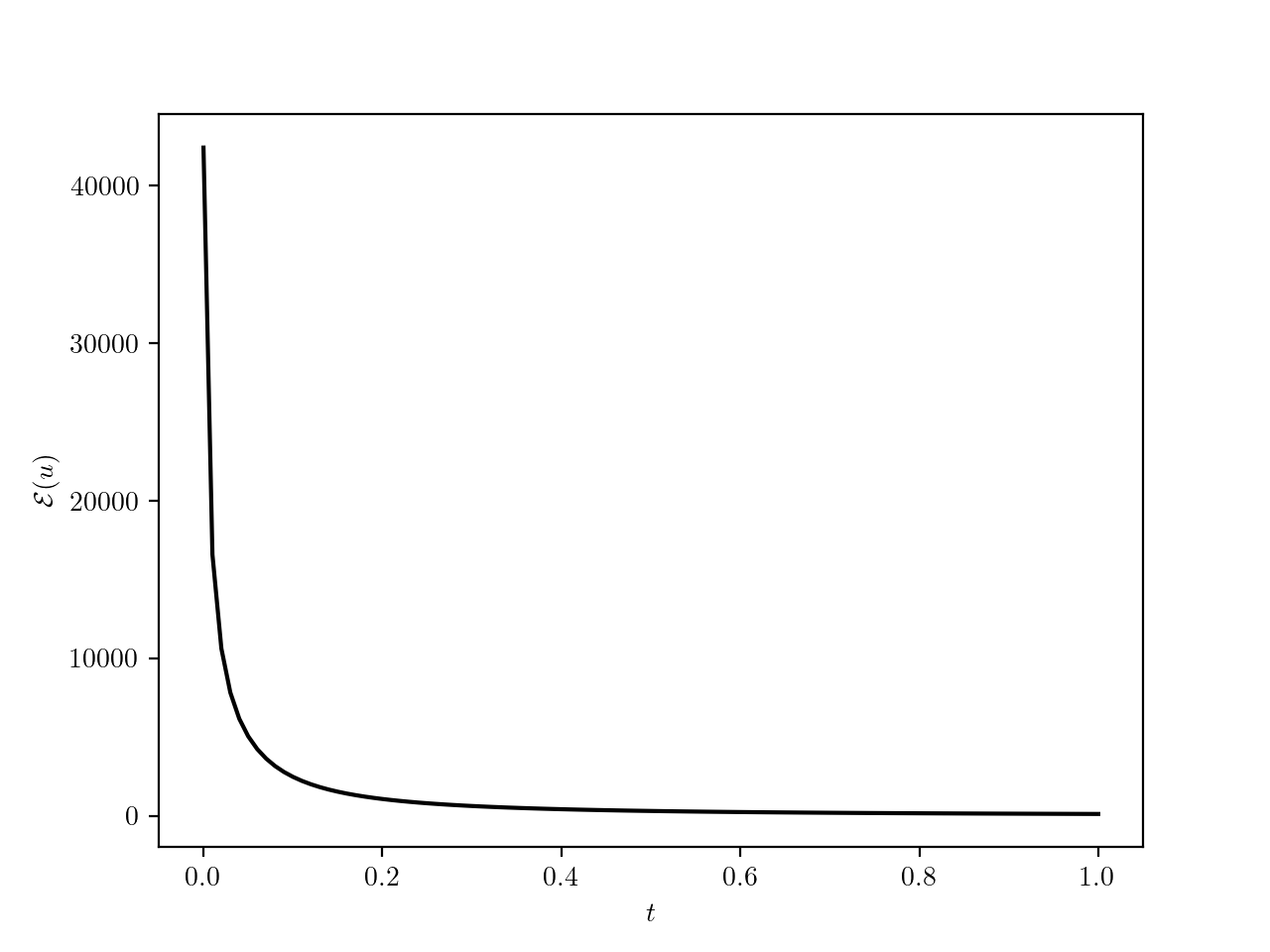}\hspace{-2ex}
\includegraphics[scale=0.35]{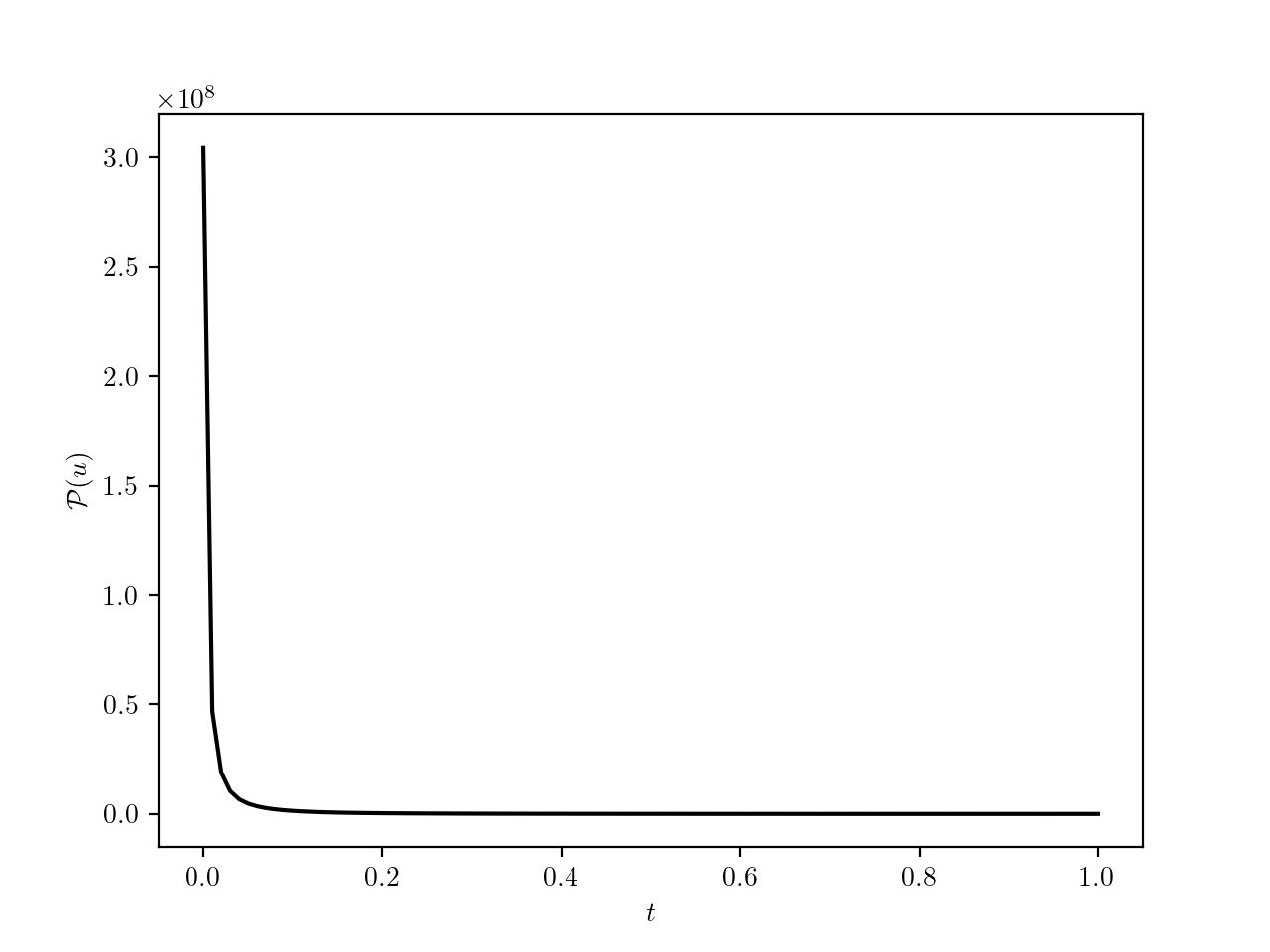}
\caption{Evolution of the Hamiltonian (left), of the enstrophy (centre), and of the palinstrophy (right).}
\label{fig:Sph_vtx}
\end{figure}
In \Cref{fig:Sph_vtx} we plot the time evolution of the Hamiltonian (on the left), enstrophy (at the centre) and palinstrophy (on the right) for this simulation. As we can see the three quantities decay monotonically, as expected.

\section{Conclusions and perspectives}\label{sec:conclusions}

We have addressed the problem of finding a geometric approximation of nonconservative evolution equations characterised by a Hamiltonian part and a gradient flow dissipation.
We have proposed a conformal variational discretisation based on the mixed formulation of the problem that ensures preservation of the geometric structure of the problem and convergence of the approximate solution to the exact one.
The focus of this work is on the spatial discretisation of problem \eqref{eq:pbm}, the investigation of numerical time integration scheme is left for future work.

To achieve sufficiently accurate solution the method might require a relatively large approximation space.
The introduction of a dynamical adaptive discretisation in
space could considerably speed up numerical simulations.
An adaptive discretisation that does not hinder the geometric structure of the problem will be subject of future investigation.

\printbibliography


\end{document}